\documentclass[12pt,twoside]{article} % leqno
\usepackage[latin1]{inputenc}
\usepackage{graphicx}
\usepackage{color}

\usepackage{enumitem}
\usepackage{latexsym}
\usepackage{amssymb}
\usepackage{epsfig}
\usepackage{hyperref}
\usepackage{amsmath}
\usepackage{natbib}
\usepackage{float}
\usepackage{dsfont}

\allowdisplaybreaks %erlaubt Seitenumbr³che in Arrays

\newcommand{\Var}{\mbox{Var}}
\def\eps{\varepsilon}
\def \R{\mathbb{R}}

\def\er{\mathbb{R}}

\def\e{\varepsilon}
\def\he{\hat{\varepsilon}}

\def\wj{\bar{w}_j}
\def\wl{\bar{w}_{\ell}}

\def\sjn{\sum_{j=1}^n}
\def\snk{\sum_{\nu=1}^kt_\nu X_{j-\nu}}

\def\beq{\begin{eqnarray*}}
\def\eeq{\end{eqnarray*}}
\def\beqn{\begin{eqnarray}}
\def\eeqn{\end{eqnarray}}

\def\hk{\hat{\kappa}_n}

\definecolor{lila}{rgb}{0.6,0.0,0.8}

\begin{document}

\title{\bf
Specification testing in nonparametric AR-ARCH models\footnote{Financial support of the DFG (Research Unit FOR 1735 {\it Structural Inference in
		Statistics: Adaptation and Effciency}) and GA\v CR 15-09663S is gratefully acknowledged.
Corresponding author: Natalie Neumeyer, University of Hamburg, Department of Mathematics, Bundesstrasse 55, 20146 Hamburg, Germany, e-mail: neumeyer@math.uni-hamburg.de}
}

\author{ Marie Hu\v{s}kovß, \\ \small Department of Statistics\\\small Charles University of Prague\!\!\!\! \and Natalie Neumeyer,\,  Tobias Niebuhr,\,  Leonie Selk
\\ \small Department of Mathematics \\\small University of Hamburg}

\maketitle

\newtheorem{theo}{Theorem}[section]
\newtheorem{lemma}[theo]{Lemma}
\newtheorem{cor}[theo]{Corollary}
\newtheorem{rem}[theo]{Remark}
\newtheorem{rems}[theo]{Remarks}
\newtheorem{prop}[theo]{Proposition}
\newtheorem{defin}[theo]{Definition}
\newtheorem{example}[theo]{Example}

\begin{abstract}
In this paper an autoregressive time series model with conditional heteroscedasticity is considered, where both conditional mean and conditional variance function are modeled nonparametrically.
A test for the model assumption of independence of innovations from past time series values is suggested. The test is based on an weighted $L^2$-distance of empirical characteristic functions. The asymptotic distribution under the null hypothesis of independence is derived and consistency against fixed alternatives is shown. A smooth autoregressive residual bootstrap procedure is suggested and its performance is shown in a simulation study.
\end{abstract}

% Running title: ...



Running title: Nonparametric AR-ARCH models

\section{Introduction}\label{introduction}
\def\theequation{1.\arabic{equation}}
\setcounter{equation}{0}

Assume we have observations  from a one-dimensional stationary weakly dependent time series $X_j$, $j\in\mathbb{Z}$.
Nonparametric modeling avoids misspecification problems and thus such models have gained much attention over the last years, see Fan and Yao (2003) and Gao (2007) for extensive overviews.
One popular possibility is to analyze data by fitting a nonparametric AR(1)-ARCH(1)-model (also called CHARN-model), i.\,e.
\[X_j=m(X_{j-1})+\sigma(X_{j-1})\e_j,\quad j\in\mathbb{Z},\]
with autoregression function $m(x)=E[X_j\mid X_{j-1}=x]$, conditional variance function $\sigma^2(x)=\mbox{Var}(X_j\mid X_{j-1}=x)$,
and innovations $\e_j$, independent from past time series values $X_{j-1},X_{j-2},\dots$. Before applying any procedure developed for a time series model like the one defined, model assumptions need to be tested.
Thus we are interested in testing the hypothesis
\[H_0:\ \e_j\mbox{ and } (X_{j-1},X_{j-2}\ldots) \mbox{ are stochastically independent.}\]
Although testing for this model assumption is essential for applications in order to obtain correct forecasts, it seems that the problem has not been considered before in the literature for the nonparametric case.
The reason is presumably that tests for hypotheses involving the innovation distribution would typically be based on the empirical distribution function of nonparametrically estimated innovations (residuals). Only recently, asymptotic results for such processes in nonparametric autoregressive models are available. M³ller et al.\ (2009) consider the above model in the homoscedastic case with constant $\sigma$. They prove an asymptotic expansion of the empirical process of residuals obtained from local-polynomial estimation of the autoregression function $m$. Further, Dette et al.\ (2009) base a test for the multiplicativity hypothesis $m=c\sigma$ on the estimated innovation distribution.
Selk and Neumeyer (2013) consider  sequential empirical process of residuals and apply it to test for a change-point in the innovation distribution.
In order to test an implication of the  null hypothesis $H_0$ one could consider, for some fixed and prespecified $k\in\mathbb{N}$,  test statistics based on an estimated  difference of the joint empirical distribution function of $\e_j$ and $ (X_{j-1},\ldots,X_{j-k}) $ and the product of the marginal distributions. Asymptotic theory could be derived similar to the considerations in M³ller et al.\ (2009), Dette et al.\ (2009), and Selk and Neumeyer (2013). Note, however,  that the assumptions for deriving asymptotic distributions of residual-based processes as in the aforementioned literature are very restrictive. To avoid unnecessarily strong assumptions we follow a different path in the paper at hand and base our test on an estimated weighted $L^2$-distance between the joint and the marginal characteristic functions of $\e_j$ and $ (X_{j-1},\ldots,X_{j-k}) $.
In an iid context a test for independence of errors and covariates in nonparametric regression models based on residual empirical characteristic functions was suggested by Hlßvka et al.\ (2011). Relatedly, in a time series context but for a parametric model  Hlßvka et al.\ (2012) test for a change in the innovation distribution of a linear autoregression model based on residual empirical characteristic functions.
Another motivation for considering the empirical characteristic functions instead of empirical distribution functions is that in other contexts it has been observed that those tests inhabit better power properties,  e.g., see Hlßvka et al.\ (2016). A survey of testing procedures based on empirical characteristic functions is given in Meintanis (2016).

The remainder of the paper is organized as follows. In section \ref{sec-null} we define our estimators and the test statistic. In section \ref{sec-asymp} we state model assumptions and give the asymptotic distribution of the test statistic under the null hypothesis, whereas consistency under fixed alternatives is discussed in section \ref{sec-alt}. A bootstrap procedure is suggested in section \ref{bootstrap}, where also the finite sample performance is investigated in a simulation study. Section \ref{conclude} concludes the paper, while all proofs are presented in an appendix.

\section{The test statistic}\label{sec-null}
\def\theequation{2.\arabic{equation}}
\setcounter{equation}{0}

Assume we have observations  $X_{-k+1},\ldots,X_n$ from the  time series $X_j$, $j\in\mathbb{Z}$, considered in section \ref{introduction}.
As test statistic for independence of innovations and past time series values we consider the weighted $L^2$-distance
\[T_n\;=\;n\int\left|\hat{\varphi}_{\he,\bar{X}_k}(t_0,t_1,\ldots,t_k)-\hat{\varphi}_{\he}(t_0)\hat{\varphi}_{\bar{X}_k}(t_1,\ldots,t_k)\right|^2W(t_0,\ldots,t_k)\,d(t_0,\ldots,t_k).\]
Here $W$ denotes some weight function fulfilling assumption \ref{A5} in Section \ref{sec-asymp}.
Furthermore
$$\hat{\varphi}_{\he,\bar{X}_k}(t_0,t_1\ldots,t_k)=\ \sjn \wj\exp\left(i\left(t_0\he_j+\snk\right)\right)$$
estimates the joint characteristic function of $\eps_j$ and $\bar{X}_{k,j}=(X_{j-1},\dots,X_{j-k})$, whereas
\begin{eqnarray*}
\hat{\varphi}_{\he}(t)&=& \sjn \wj \exp\left(it\he_j\right),\\
\hat{\varphi}_{\bar{X}_k}(t_1,\ldots,t_k)&=&\sum_{j=1}^n\wj\exp\left(i\sum_{\nu=1}^kt_{\nu}X_{j-\nu}\right)
\end{eqnarray*}
estimate the marginal characteristic functions of $\e_j$ and $\bar{X}_{k,j}$, respectively. Here the weights
 are defined as
$\wj=w_n(X_{j-1})/(\sum_{l=1}^nw_n(X_{l-1}))$, where we choose a weight function $w_n(x)=I_{[-a_n,a_n]}(x)$ for some sequence $a_n\to \infty$. Here and throughout $I_A$ denotes the indicator function of set $A$. Other weight functions $w_n:\er\to[0,1]$ which vanish outside $[-a_n,a_n]$ are possible as well but require slightly adapted assumptions.  The weights are included in the definition of the empirical characteristic functions to avoid problems of kernel estimation in areas where only few data are available.
Furthermore the residuals are defined as
$\hat{\e}_j=(X_j-\hat{m}(X_{j-1}))/\hat{\sigma}(X_{j-1})$
and we use Nadaraya-Watson type estimators for the conditional mean and variance functions,
\begin{eqnarray*}
\hat{m}(x)&=&\frac{\frac 1{nc_n}\sum_{j=1}^nK(\frac{x-X_{j-1}}{c_n})X_j}{\hat f_X(x)}\\
\hat{\sigma}^2(x)&=&\frac{\frac 1{nc_n}\sum_{j=1}^nK(\frac{x-X_{j-1}}{c_n})(X_j-\hat m(x))^2}{\hat f_X(x)}
\end{eqnarray*}
with kernel function $K$ and sequence of bandwidths $c_n, n\in\mathbb{N}$. Here
$$\hat f_X(x)=\frac 1{nc_n}\sum_{j=1}^nK\Big(\frac{x-X_{j-1}}{c_n}\Big)$$
denotes a kernel estimator for the marginal density $f_X$ of $X_j$. See, e.\,g.,  Robinson (1983), Masry and Tj$\o$stheim (1995), Hõrdle and Tsybakov (1997) and Hansen (2008)
for properties of these estimators in the time series context.

\section{Assumptions and asymptotic results under the null hypothesis}\label{sec-asymp}
\def\theequation{3.\arabic{equation}}
\setcounter{equation}{0}

Under the null hypothesis we state the following assumptions. Please note that throughout we write $\bold{t}=(t_0,t_1,\dots,t_k)$ and use the notation $g(\bold{t})$ for simplicity also for functions $g$ that only depend on $(t_1,\dots,t_k)$ (see e.\,g.\ $\psi(\bold{t},x)$ from assumption \ref{Apsi}).
\begin{enumerate}[label=(\textbf{A\arabic{*}})]
\item \label{A1} The process $(X_t)_{t\in\mathbb{Z}}$ is strictly stationary and $\alpha$-mixing with mixing coefficient $\alpha$ that satisfies  $\alpha(i)\leq Ai^{-\beta}$ for some $A<\infty$ and $\beta>\frac{1+ (s-1)(2+1/q)}{s-2}$ for some $q>0$, where $s>2$ and  $E|X_0|^s<\infty$.
\\
$X_1$ has bounded marginal density $f_X$ such that for some constant $B_1$, $$\sup_x E(|X_1|^s| X_0=x) f_X(x) \leq B_1.$$
 Furthermore $(X_0,X_j)$ has bounded joint density $f_j$ and there exists a constant $B_2$, such that for  some $j^*$,
$$
\sup_{x_0,x_j} E(|X_1 X_{j+1}|X_0=x_0, X_j=x_j) f_j(x_0, x_j) \leq B_2
$$
for  all $j\geq j^*$.
\item \label{Am}
Let $m$, $\sigma^2$ and $f_X$ be differentiable. Let there  exist some $r\in (0,\infty)$ such that  the functions $m,m',\sigma^2, (\sigma^2)',\frac{1}{\sigma^2},\frac{1}{f_X}$ and $f_X'$ are of order $O((\log n)^r)$ uniformly on the interval $I_n=[-a_n-Cc_n,a_n+Cc_n]$ (with $C$ from assumption \ref{A3}).
Further we assume Lipschitz continuity of the derivatives $f_X'$, $m'$ and $(\sigma^2)'$ in the following sense, 
$$\sup_{x,y\in I_n\atop |x-y|\leq c_n}|g(x)-g(y)|=O(c_n(\log n)^r) \mbox{ for } g\in \{f_X',m',(\sigma^2)'\}.$$
\item \label{A2} The innovations $(\varepsilon_t)_{t\in\mathbb{Z}}$ are independent, centered and identically distributed. For each $t\in\mathbb{Z}$,  $\varepsilon_t$ is independent from the past $X_{t-1},X_{t-2},\dots$.

For some $\delta>\frac{2}{\beta-2}$ let $E[|\varepsilon_1|^{(2+2\delta)\vee 4}]<\infty$ and
$\sup_{x\in I_n} E[|\varepsilon_j|^{2(1+\delta)}\mid X_0=x]=O((\log n)^r)$ uniformly in $j$ with $r$ and $I_n$ from assumption \ref{Am}.

\item \label{Apsi} Define $\psi(\bold{t},x)=E[Y_1(\bold{t})|X_0=x]-E[Y_1(\bold{t})]$ with $Y_1(\bold{t})=\cos(\sum_{\nu=1}^k t_\nu X_{1-\nu})$ and 
assume that 
$$\displaystyle \sup_{x,z\in I_n\atop |x-z|\leq Cc_n}\sup_\bold{t} |\psi(\bold{t},x)-\psi(\bold{t},z)|=O((\log n)^rc_n^d)$$ for some $d>0$ with $r$ from assumption \ref{Am} and $C$ from assumption \ref{A3}.
Assume the same condition holds for 
$\tilde\psi(\bold{t},x)=E[Z_1(\bold{t})|X_0=x]-E[Z_1(\bold{t})]$ with $Z_1(\bold{t})=\sin(\sum_{\nu=1}^k t_\nu X_{1-\nu})$.

\item \label{A3} The kernel $K$ is a symmetric and Lipschitz continuous density with compact support $[-C,C]$ and $\int K(u)u\,du=0$.

\item\label{A4}
For $q$, $s$ and $\beta$ from \ref{A1} we have $a_n=O(n^{1/(2q)}\log n)$, and for $\theta=\frac{\beta-2-\frac 1q-\frac{1+\beta}{s-1}}{\beta+2-\frac{1+\beta}{s-1}}$ it holds that $\log n=o(n^\theta c_n)$.
Let
$$a_n^*=\left(\frac{\log n}{nc_n}\right)^{1/2}+c_n^2,$$
then $a_n^*=O(\Delta_nn^{-1/4})$ with $\Delta_n=\inf_{|x|\leq a_n}f_X(x)$.

\item\label{Aneu}
 Let the sequence of bandwidths fulfill $nc_n^{2}(\log n)^{-D}\to\infty$, $nc_n^{4}(\log n)^{D}\to0$ for all $D>0$.
 
\item\label{A5} The weight function $W$ is nonnegative and symmetric such that $W(\pm t_0,\pm t_1,\dots,\pm t_k)$\linebreak$=W(t_0,\dots,t_k)$. Further $\int t_0^4W(t_0,\ldots,t_k)d(t_0\ldots,t_k)<\infty$.

\end{enumerate}

\begin{rem}
Apart from the typical assumptions on the kernel, bandwidths and weight functions we need smoothness assumptions on the unkown functions as well as moment assumptions and the mixing property, e.\,g.\ in order to obtain uniform rates of convergence for the kernel estimators, similar to Hansen (2008). 
Note that for \ref{A4} and \ref{Aneu} both to be satisfied one needs $\theta >\frac 14$. 
\end{rem}

We have the following asymptotic distribution of the test statistic under the null.

\begin{theo}\label{theo1} Under the assumptions \ref{A1}--\ref{A5} the test statistic $T_n$ converges in distribution to $T=\int_{\R^{k+1}} S^2(\bold{t}) W(\bold{t})\, d\bold{t}$, where $S(\bold{t})$, $\bold{t}\in\R^{k+1}$, denotes a centered Gaussian process with the same covariance structure as
\begin{eqnarray*}
&&\tilde{S}(t_0,\dots,t_k)\\
&=&\Big(\cos(t_0\eps_1)-E\big[\cos(t_0\eps_1)\big]\Big)\Big(Y_1(\bold{t})+Z_1(\bold{t})-E[Y_1(\bold{t})+Z_1(\bold{t})]\Big)\\
&& {} +\Big(\sin(t_0\eps_1)-E[\sin(t_0\eps_1)]\Big)\Big(Y_1(\bold{t})-Z_1(\bold{t})-E[Y_1(\bold{t})-Z_1(\bold{t})]\Big)
\\
&& {}+t_0\Big(\eps_1 E[\sin(t_0\eps_1)]+\frac{1}{2}(\eps_1^2-1)E[\sin(t_0\eps_1)\eps_1]\Big)\Big(E[Y_1(\bold{t})+Z_1(\bold{t})|X_{0}]-E[Y_1(\bold{t})+Z_1(\bold{t})]\Big)\\
&& {} -t_0\Big(\eps_1 E[\cos(t_0\eps_1)]+\frac{1}{2}(\eps_1^2-1)E[\cos(t_0\eps_1)\eps_1]\Big)\Big(E[Y_1(\bold{t})-Z_1(\bold{t})|X_{0}]-E[Y_1(\bold{t})-Z_1(\bold{t})]\Big).
\end{eqnarray*}
%where $Y_1(\bold{t})=\cos(\sum_{\nu=1}^k t_\nu X_{1-\nu})$ and $Z_1(\bold{t})=\sin(\sum_{\nu=1}^k t_\nu X_{1-\nu})$.
\end{theo}

\medskip

The proof is given in the appendix. An asymptotic level-$\alpha$ test is obtained by rejecting $H_0$ whenever $T_n>c_{1-\alpha}$, where $P(T>c_{1-\alpha})=\alpha$.
Due to the complicated distribution of $T$ we suggest a bootstrap procedure to estimate the critical value $c_{1-\alpha}$ in section \ref{bootstrap}.

\begin{rems}
{\bf (a)} The replacement of true but unknown innovations $\e_j$ by the estimated residuals $\hat\e_j$ changes the asymptotic distribution drastically. Were the true innovations  known and used in the test statistic instead of residuals the statistic $\tilde S$ in Theorem \ref{theo1} would simplify to
\begin{eqnarray*}
\tilde{S}(t_0,\dots,t_k)
&=&\Big(\cos(t_0\eps_1)-E\big[\cos(t_0\eps_1)\big]\Big)\Big(Y_1(\bold{t})+Z_1(\bold{t})-E[Y_1(\bold{t})+Z_1(\bold{t})]\Big)\\
&& {} +\Big(\sin(t_0\eps_1)-E[\sin(t_0\eps_1)]\Big)\Big(Y_1(\bold{t})-Z_1(\bold{t})-E[Y_1(\bold{t})-Z_1(\bold{t})]\Big)
.\end{eqnarray*}
{\bf (b)} If the aim is to test for independence of innovations and past time series values in a (homoscedastic) AR(1) model
$$X_j=m(X_{j-1})+\e_j,$$
one simply sets $\hat\sigma\equiv 1$ in the definition of the residuals. Then the statistic $\tilde S$ in Theorem \ref{theo1} changes to
\begin{eqnarray*}
\tilde{S}(t_0,\dots,t_k)
&=&\Big(\cos(t_0\eps_1)-E\big[\cos(t_0\eps_1)\big]\Big)\Big(Y_1(\bold{t})+Z_1(\bold{t})-E[Y_1(\bold{t})+Z_1(\bold{t})]\Big)\\
&& {} +\Big(\sin(t_0\eps_1)-E[\sin(t_0\eps_1)]\Big)\Big(Y_1(\bold{t})-Z_1(\bold{t})-E[Y_1(\bold{t})-Z_1(\bold{t})]\Big)\\
&& {}+t_0\eps_1 E[\sin(t_0\eps_1)]\Big(E[Y_1(\bold{t})+Z_1(\bold{t})|X_{0}]-E[Y_1(\bold{t})+Z_1(\bold{t})]\Big)\\
&& {} -t_0\eps_1 E[\cos(t_0\eps_1)]\Big(E[Y_1(\bold{t})-Z_1(\bold{t})|X_{0}]-E[Y_1(\bold{t})-Z_1(\bold{t})]\Big)
.\end{eqnarray*}
{\bf (c)}  As mentioned in the introduction alternative testing procedures would be given by, e.\,g., Kolmogorov-Smirnov or CramÚr-von Mises type statistics based on the $\hat F_{\he,\bar{X}_k}-F_{\he}\otimes F_{\bar{X}_k}$, i.\,e.\ the weighted empirical joint distribution function of $\hat\e_j$ and $\bar{X}_{k,j}=(X_{j-1},\dots,X_{j-k})$ ($j=1,\dots,n$) and the product of the marginals.
Following M³ller et al.\ (2009), Dette et al.\ (2009), and Selk and Neumeyer (2013) to derive the asymptotic distribution would, however, require stronger assumptions on the data generating process.
\end{rems}

\section{Fixed alternatives}\label{sec-alt}
\def\theequation{4.\arabic{equation}}
\setcounter{equation}{0}

Note that by construction the test statistic $T_n$ cannot detect alternatives where the innovation $\e_j$ is independent of $ (X_{j-1},\ldots,X_{j-k})$, but depends on some $X_{j-\ell}$ for $\ell >k$.
However, the test is consistent against any fixed alternative
$$H_1:\e_j\mbox{ and } X_{j-\ell}\mbox{ are stochastically dependent for some }\ell\in\{1,\dots,k\} $$
under the following model.
Assume that $(X_j)_{j\in \mathbb{Z}}$ is a strictly stationary and weakly dependent time series that fulfills assumption \ref{A1}. Further define $m(x)=E[X_{j+1}\mid X_j=x]$ and $\sigma^2(x)=\Var(X_{j+1}\mid X_j=x)$. Let $m$, $\sigma^2$ and the marginal density $f_X$ fulfill assumption \ref{Am}. Let the kernel, weight function and sequence of bandwidth fulfill
\ref{A3}--\ref{A5}. Then we have the following result.

\begin{theo}\label{theo2} Under the assumptions listed in this section, $T_n/n$  converges to
\[\tilde T=\int\left|{\varphi}_{\e,\bar{X}_k}(t_0,t_1,\ldots,t_k)-{\varphi}_{\e}(t_0){\varphi}_{\bar{X}_k}(t_1,\ldots,t_k)\right|^2W(t_0,\ldots,t_k)\,d(t_0,\ldots,t_k)\]
in probability, where ${\varphi}_{\e,\bar{X}_k}$ is the joint characteristic function of $\e_j$ and $(X_{j-1},\dots,X_{j-k})$, and $\varphi_{\e}$ and $\varphi_{\bar{X}_k}$ are the corresponding marginal characteristic functions.
\end{theo}

The proof is given in the appendix.
Note that under $H_1$ one has $\tilde T>0$ and hence $T_n\longrightarrow \infty$ for $n\to\infty$.

\bigskip

From rejection of $H_0$ one should conclude that the AR(1)-ARCH(1) model is not suitable to describe the data. Possible reasons are explained in the following example.

\begin{example}
{\bf (a)}
Consider the conditional distribution of $\e_j$, given $X_{j-1}$. The first two moments of this distribution do not depend on $X_{j-1}$ by construction. Higher order moments could depend on $X_{j-1}$, i.\,e.\ $E[\e_j^\ell\mid X_{j-1}]=h_\ell (X_{j-1})$ for some $\ell\geq 3$.
In the simulation study we will consider a skew normal innovation distribution with mean zero, variance one and skewness dependent on $X_{j-1}$.

{\bf (b)} The conditional distribution of $\e_j$, given $\bar{X}_k=(X_{j-1},\dots,X_{j-k})$ may still depend on $\bar{X}_k$.
If this distribution does still depend on the first component $X_{j-1}$, but only on this component, modeling the autoregression and conditional variance function with lag 1 is appropriate, but one should not apply any procedures that assume independence of innovations and past time series values.

{\bf (c)} An AR($\ell$)-ARCH($\ell$) model could be appropriate for the data for some $\ell>1$, i.\,e.\
$$X_j=\tilde m(X_{j-1},\dots,X_{j-\ell})+\tilde \sigma(X_{j-1},\dots,X_{j-\ell})\eta_j$$
with innovations $\eta_j$ independent from $X_{j-1},X_{j-2},\dots$.
\end{example}

\section{Bootstrap and finite sample performance}\label{bootstrap}
\def\theequation{5.\arabic{equation}}
\setcounter{equation}{0}

%%%%
In this section we investigate the finite-sample performance of our test by simulations.
Due to the complicated limiting distribution of $T$ from Theorem \ref{theo1}, we suggest to use a smooth autoregressive residual bootstrap instead.
%%%%%
Our bootstrap strategy is as follows.\\
Firstly, based on the estimators as introduced in section \ref{sec-null},
generate bootstrap innovations $\e_j^*$ from a smooth estimate of the innovation distribution, i.\,e. given the original data $X_{-k+1},\ldots,X_n$ the distribution of $\e_j^*$ reads
$$F_n(x)=\frac{1}{n}\sum_{i=1}^n L\Big(\frac{x-\tilde\e_i}{h_n}\Big)$$
where $h_n$ denotes a positive bandwidth, $L$ is some smooth distribution function and $\tilde \e_1,\dots,\tilde \e_n$ denote the standardized versions of the residuals $\hat\e_1,\dots,\hat\e_n$.
Secondly, compute the bootstrap process via
$$X_j^*=\hat m(X_{j-1}^*)+\hat \sigma(X_{j-1}^*)\e_j^*,\;\;j=1,\ldots,n,$$
with some starting value $X_0^*$ and a sufficiently large number of forerunnings to ensure the process is in balance.
Thirdly, calculate the bootstrap analogue of the test statistic $T_n$, say $T_n^*$.\\
Frequent repetitions of these steps give the distribution of $T_n^*$ which approximates the distribution of $T_n$.
By using the $1-\alpha$ percentile of the distribution of $T_n^*$, say $c_{1-\alpha}^*$, 
the hypothesis of independence then is rejected if $T_n>c^*_{1-\alpha}$.
It is worth noting that, given the original data $X_{-k+1},\ldots,X_n$,  the bootstrap innovations $\e_j^*$ are independent of $X_{j-1}^*,X_{j-2}^*,\dots$ and thus the bootstrap data fulfills the null hypothesis.\\
%\textcolor{blue}{Furthermore note that drawing $\e_j^*$ with replacement from (standardized) residuals is not appropriate as then the bootstrap process would not be mixing (REF!).}\\
The simulations are restricted to the hypothesis '$H_0$: $\varepsilon_j$ and $X_{j-1}$ are stochastically independent', i.e. only the case $k=1$ is investigated.
To examine the performance of the test for finite sample sizes, we consider the following two AR-ARCH models:
\begin{align*}
(i)\;m(x)=0.9x,\; \sigma(x)\equiv 1,\;\;\; (ii)\;m(x)\equiv 0,\; \sigma(x)=\sqrt{1+0.25x^2}.
\end{align*}
Obviously, model (i) corresponds to an AR series and model (ii) represents an ARCH model.
For both models, the performance under the null and under the alternative is investigated.
We distinguish the null from the alternative by the choice of the innovation sequence. Under the null we use standard normally distributed innovations. Under the alternative we choose standardized skew normally distributed innovations, where the skewness parameter depends on past time series values. In particular the skewness parameter of $\varepsilon_{t+1}$ was set to $10X_t^2$ for all relevant time points $t$, according to the notation of Fern\'andez and Steel (1998).

\scriptsize{\textbf{Table 1:} Rejection probabilities for the AR model (i) under the null hypothesis (left) and the under the alternative (right).}
\begin{minipage}{0.5\textwidth}
\begin{align*}
\begin{array}{lrrr}
\hline
&\alpha=0.01&\alpha=0.05&\alpha=0.1\\
n= 50&0.0000&0.0275&0.0475\\
n=100&0.0000&0.0125&0.0225\\
n=200&0.0005&0.0175&0.0425\\
n=300&0.0000&0.0200&0.0375\\
n=400&0.0000&0.0125&0.0425\\
\hline
\end{array}
\end{align*}
\end{minipage}
\begin{minipage}{0.5\textwidth}
\begin{align*}
\begin{array}{lrrr}
\hline
&\alpha=0.01&\alpha=0.05&\alpha=0.1\\
n= 50&0.0150&0.0600&0.1125\\
n=100&0.0450&0.1375&0.2200\\
n=200&0.1550&0.3225&0.5200\\
n=300&0.5100&0.7700&0.8675\\
n=400&0.8350&0.9525&0.9900\\
\hline
\end{array}
\end{align*}
\end{minipage}\\
\normalsize{}

\scriptsize{\textbf{Table 2:} Rejection probabilities for the ARCH model (ii) under the null hypothesis (left) and under the alternative (right).}
\begin{minipage}{0.5\textwidth}
\begin{align*}
\begin{array}{lrrr}
\hline
&\alpha=0.01&\alpha=0.05&\alpha=0.1\\
n= 50&0.0075&0.0150&0.0225\\
n=100&0.0000&0.0125&0.0325\\
n=200&0.0150&0.0325&0.0700\\
n=300&0.0025&0.0225&0.0625\\
n=400&0.0075&0.0475&0.0825\\
\hline
\end{array}
\end{align*}
\end{minipage}
\begin{minipage}{0.5\textwidth}
\begin{align*}
\begin{array}{lrrr}
\hline
&\alpha=0.01&\alpha=0.05&\alpha=0.1\\
n= 50&0.0325&0.0900&0.1525\\
n=100&0.1000&0.2125&0.2775\\
n=200&0.3125&0.4750&0.5275\\
n=300&0.4900&0.5750&0.6375\\
n=400&0.5575&0.6025&0.6450\\
\hline
\end{array}
\end{align*}
\end{minipage}\\
\hfill\\
\normalsize{}
%%%%%%%%%%%%%%%%%%%%%%%%%%%%%%%%%%%%%%%%%%%%%%%
Tables 1 and 2 state the rejection probabilities for 400 Monte Carlo simulations each with 400 bootstrap repetitions for several sample sizes $n$ and significance levels $\alpha$.
We chose $L$ as the standard normal distribution, $h_n$ was set to $n^{-1/4}$ for reasons given in Neumeyer (2006), and the bandwidth $c_n$ was chosen by Silverman's rule by thumb, see Silverman (1986).\\
The tables show that under the null hypothesis the test yields the given level of significance.
While for model (ii) the test performance is very likely, for model (i) the test seems to be somehow over-conservative for the sample sizes used. 
Under the alternative the test power increases with increasing sample size in both models as to be expected.
It is worth to note that the test power increases faster for model (i) than for model (ii).
Altogether, the procedure performs satisfying in our simulations, however, it has to be noticed that the test performance depends on  the time series at hand.
% and whether the null hypothesis is fulfilled or not.
\\
For practitioners the computation of the test statistic $T_n$, and $T_n^*$ respectively, might be challenging.
For that reason, we suppose using another representation of $T_n$, and $T_n^*$, which avoids for solving complicated integrals.
The alternative representation is stated in the following lemma.

\begin{lemma}\label{fourier}
Let $\mathcal F[V](x)$ denote the Fourier transformation of $V$ at point $x$.
Under Assumption $\ref{A5}$ and if $W$ yields $W(t_0,\ldots,t_k)=V_0(t_0)\prod_{i=1}^kV_i(t_i)$,
it holds $\mathcal F[V](x)=\int\cos(tx)V(t)dt$ and
the test statistic $T_n$ can be represented by
\begin{align*}
T_n=\label{fourierTn}
&
n\sum_{s_1,s_2=1}^n%\sum_{s_2=1}^n\sum_{s_3=1}^n\sum_{s_4=1}^n
\bar w_{s_1}\bar w_{s_2}%\bar w_{s_3}\bar w_{s_4}
\mathcal F[V_0](\hat\varepsilon_{s_1}-\hat\varepsilon_{s_2})
\sum_{s_3,s_4=1}^n
\bar w_{s_3}\bar w_{s_4}
\prod_{j=1}^k\mathcal F[V_j](X_{s_3-j}-X_{s_4-j})\\
&+
n\sum_{s_1,s_2=1}^n%\sum_{s_2=1}^n
\bar w_{s_1}\bar w_{s_2}
\mathcal F[V_0](\hat\varepsilon_{s_1}-\hat\varepsilon_{s_2})
\prod_{j=1}^k\mathcal F[V_j](X_{s_1-j}-X_{s_2-j})
\notag\\
&-
2n\sum_{s_1,s_2,s_3=1}^n%\sum_{s_2=1}^n\sum_{s_3=1}^n
\bar w_{s_1}\bar w_{s_2}\bar w_{s_3}
\mathcal F[V_0](\hat\varepsilon_{s_1}-\hat\varepsilon_{s_2})
\prod_{j=1}^k\mathcal F[V_j](X_{s_1-j}-X_{s_3-j}).
\notag
\end{align*}
\end{lemma}
Since the choice of the weighting function $W$ belongs to the user, the additional assumption on its multiplicative form is very weak.
If one further chooses $W$ such that the Fourier transformations of the corresponding functions $V_i$, $i=0,\ldots,k$, are known, the test statistic $T_n$  can straightforwardly be computed.
Even more important, the implementation then simplifies a lot since the computation of the $(k+1)$-fold integral is omitted.
\begin{example}
Some choices of $W$ fulfilling the assumptions of the lemma are:
\begin{itemize}
\item[\bf{(a)}] $W(t_0,\ldots,t_k)%=e^{-\sum_{j=0}^{k}\gamma_j|t_j|}
=e^{-\gamma_0 |t_0|}\prod_{j=1}^ke^{-\gamma_j|t_j|}$, where the Fourier transformation of $V_j(t_j):=e^{-\gamma_j|t_j|}$, $j=0,\ldots,k$, is given by $\mathcal F[V_j](t_j)=\frac{2\gamma_j}{\gamma_j^2+4\pi^2x^2}$.
\item[\bf{(b)}] $W(t_0,\ldots,t_k)%=e^{-\sum_{j=0}^{k}\gamma_j t_j^2}
=e^{-\gamma_0 t_0^2}\prod_{j=1}^ke^{-\gamma_j t_j^2}$, where the Fourier transformation of $V_j:=e^{-\gamma_j t_j^2}$, $j=0,\ldots,k$, is given by $\mathcal F[V_j](x)=\sqrt{\frac{\pi}{\gamma_j}}e^{-(\pi x)^2/\gamma_j}$.
\end{itemize}
\end{example}
%%%%%%%%%%%%%%%%%%%%%%%%%%%%%%%%%%%%%%%%%%%%%%%%%%%%%%%

\section{Concluding remarks and outlook}\label{conclude}
\def\theequation{7.\arabic{equation}}
\setcounter{equation}{0}

In this paper we suggested a test for independence of innovations and past time series observations in an AR-ARCH model, where both the conditional mean and conditional volatility function are modeled nonparametrically. The test is based on empirical characteristic functions.
For simplicity of presentation we considered the AR(1)-ARCH(1) case. However, generalizations to AR($p$)-ARCH($p$) models are straightforward, while then local polynomial estimators for the mean and variance function should be used. Facing the curse of dimensionality also semiparametric models might be of interest, see e.\,g.\ Yang et al.\ (1999) for a model with an additive autoregression function and multiplicative volatility function.
Including covariates is possible as well. Then one considers a model of type
$X_j=m(T_j)+\sigma(T_j)\eps_j,$
where the vector $T_j$ may include past observations. Testing independence of $\eps_j$ from $T_j,T_{j-1},\dots$ would be of interest here and can be conducted in an analogous manner.

A question related to the one considered in the paper at hand is whether the innovations really form an iid sequence. Corresponding tests for parametric times series models
have been considered by Ghoudi et al.\ (2001), among others.
Presumably with the methods developed in the paper at hand, such hypotheses tests for nonparametric time series models can be derived. We leave the consideration for future research.

\section*{References}

\begin{description}
\item Dette, H., Pardo-Fernßndez, J. C. \& Van Keilegom, I. (2009). {\em Goodness-of-Fit Tests for Multiplicative Models with Dependent Data}. Scand. J. Statist. 36, 782--799.

\item Fan, J. \& Yao, Q. (2003). {\em Nonlinear Time Series: Nonparametric and Parametric Methods}.
Springer Series in Statistics, New York.

\item Fern\'andez, C. and Steel, F.J. (1998). \textit{On Bayesian modeling of fat tails and skewness}. J. American Stat. Assoc. 93, 259--371.

\item Gao, J. (2007). {\em Nonlinear Time Series: Semiparametric and Nonparametric Methods}. Chapman \& Hall/CRC, Boca Raton.

\item Ghoudi, K., Kulperger, R. J. \& RÚmillard, B. (2001). {\em A Nonparametric Test of Serial Independence for Time Series and Residuals.} J. Multivariate Anal. 79, 191--218.

 \item Hansen, B.E. (2008). \textit{Uniform Convergence Rates for Kernel Estimation with Dependent Data}. Econom. Theory 24, 726-748.

\item Hõrdle, W. \& Tsybakov, A. (1997). \textit{Local polynomial estimators of the volatility function in nonparametric autoregression.} J. Econometrics 81, 223-242.

\item Hlßvka, Z., Hu\v{s}kovß, M., Kirch, C. \& Meintanis, S. G. (2012). {\em Monitoring changes in the error distribution of autoregressive models based on Fourier methods.} TEST 21, 605--634.

\item Hlßvka, Z., Hu\v{s}kovß, M., Kirch, C. \& Meintanis, S. G. (2014). {\em Fourier-Type Tests Involving Margingale Difference Processes.} Econometric Reviews.\\ DOI:
10.1080/07474938.2014.977074

\item Hlßvka, Z., Hu\v{s}kovß, M., Kirch, C. \& Meintanis, S. G. (2016). {\em Bootstrap procedures for on-line monitoring of changes
in autoregressive models.}  Commun. Stat. Simul. Comput. 45, 2471--2490.

\item Hlßvka, Z., Hu\v{s}kovß, M. \& Meintanis, S. G. (2011). {\em Tests for independence in non-parametric heteroscedastic regression models.} J. Multivariate Anal. 102, 816--827.

\item Ibragimov, I.A. \& Khasminskii, R.Z. (1981). {\em Statistical Estimation: Asymptotic Theory.} Springer, New York.

\item Masry, E. \& Tj$\o$stheim, D. (1995). \textit{Nonparametric estimation and identification of nonlinear ARCH time series.} Econometric Theory 11, 258--289.

\item Meintanis, S. G. (2016). \textit{A review of testing procedures based on the empirical characteristic function.} South African Statistical Journal 50, 1--14.

\item M³ller, U. U., Schick, A. \& Wefelmeyer, W. (2009). {\em Estimating the innovation distribution in nonparametric autoregression}. Probab. Theory Relat. Fields 144, 53--77.

\item Neumeyer, N. (2009). {\em Smooth residual bootstrap for empirical processes of nonparametric regression residuals.} Scand. J. Statist. 36,  204--228.

\item Robinson, P. M. (1983). \textit{Nonparametric estimators for time series.} J. Time Ser. Anal. 4, 185-207.

\item Selk, L. \& Neumeyer, N. (2013). \textit{Testing for a change of the innovation distribution in nonparametric autoregression - the sequential empirical process approach}. Scand. J. Statist. 40, 770--788

\item Silverman, B.W. (1986). {\em Density estimation for statistics and data analysis.} Chapman and Hall, New York.

\item Su, L. \& Xiao, Z. (2008). \textit{Testing structural change in time-series nonparametric regression models.} Statistics and Its Interface. Vol. 1, 347-366.

\item Sun, S. \& Chiang, C-Y. (1997). \textit{Limiting behavior of the perturbed empirical distribution functions evaluated at U-statistics for strongly mixing sequences of random variables.} J. Appl. Math. Stoch. Anal.  10, 3--20.

\item Yang, L., Hõrdle, W. \& Nielsen, J. P. (1999). {\em Nonparametric autoregression with multiplicative volatility and additive mean.} J. Time Ser. Anal. 20, 579--604.

\item Yokoyama, R. (1980). \textit{Moment Bounds for Stationary Mixing Sequences.}
Z. Wahrscheinlichkeitstheorie verw. Gebiete 52, 45--57.

\end{description}

\begin{appendix}

\section{Proofs: main results}
\def\theequation{A.\arabic{equation}}
\setcounter{equation}{0}

Throughout the proof $D$ denotes some generic positive constant, independent of $\bold{t}$, that may differ from line to line. 

\bigskip

{\bf Proof of Theorem \ref{theo1}. }

Note that for the test statistic we have
\begin{eqnarray*}
T_n &=& n\int \Big|\sum_{j=1}^n\wj   \exp\Big(i\Big(t_0\hat\eps_j+\sum_{\nu=1}^k t_\nu X_{j-\nu}\Big)\Big)
-\sum_{j=1}^n\sum_{\ell=1}^n\wj\wl \exp\Big(i\Big(t_0\hat\eps_j+\sum_{\nu=1}^k t_\nu X_{\ell-\nu}\Big)\Big)\Big|^2\\
&&{}\qquad\times W(t_0,\dots,t_k)\, d(t_0,\dots,t_k)\\
&=& n\int \Big[\Big(\sum_{j=1}^n\wj  \cos\Big(t_0\hat\eps_j+\sum_{\nu=1}^k t_\nu X_{j-\nu}\Big)-\sum_{j=1}^n\sum_{\ell=1}^n\wj\wl\cos\Big(t_0\hat\eps_j+\sum_{\nu=1}^k t_\nu X_{\ell-\nu}\Big) \Big)^2\\
&&{}\qquad
+\Big(\sum_{j=1}^n\wj  \sin\Big(t_0\hat\eps_j+\sum_{\nu=1}^k t_\nu X_{j-\nu}\Big)-\sum_{j=1}^n\sum_{\ell=1}^n\wj\wl\sin\Big(t_0\hat\eps_j+\sum_{\nu=1}^k t_\nu X_{\ell-\nu}\Big)\Big)^2
\Big]\\
&&{}\qquad\times W(t_0,\dots,t_k)\, d(t_0,\dots,t_k)
\end{eqnarray*}
and with the addition theorems for trigonometric functions one obtains
\begin{eqnarray*}
T_n &=&
n\int \Bigg\{\Big[\sum_{j=1}^n\wj\cos(t_0\hat\eps_j)\Big(\cos\Big(\sum_{\nu=1}^k t_\nu X_{j-\nu}\Big)-\sum_{\ell=1}^n\wl\cos\Big(\sum_{\nu=1}^k t_\nu X_{\ell-\nu}\Big)\Big)\\
&&\qquad
- \sum_{j=1}^n\wj\sin(t_0\hat\eps_j)\Big(\sin\Big(\sum_{\nu=1}^k t_\nu X_{j-\nu}\Big)-\sum_{\ell=1}^n\wl\sin\Big(\sum_{\nu=1}^k t_\nu X_{\ell-\nu}\Big)\Big)
\Big]^2\\
&&{}\qquad +  \Big[\sum_{j=1}^n\wj\sin(t_0\hat\eps_j)\Big(\cos\Big(\sum_{\nu=1}^k t_\nu X_{j-\nu}\Big)-\sum_{\ell=1}^n\wl\cos\Big(\sum_{\nu=1}^k t_\nu X_{\ell-\nu}\Big)\Big)\\
&&{}\qquad +  \sum_{j=1}^n\wj\cos(t_0\hat\eps_j)\Big(\sin\Big(\sum_{\nu=1}^k t_\nu X_{j-\nu}\Big)-\sum_{\ell=1}^n\wl\sin\Big(\sum_{\nu=1}^k t_\nu X_{\ell-\nu}\Big)\Big)\Big]^2\Bigg\}\\
&&{}\quad\times W(t_0,\dots,t_k)\, d(t_0,\dots,t_k).
\end{eqnarray*}
From assumption \ref{A5} by symmetry properties of cosine and sine we obtain
\begin{eqnarray*}
T_n &=&
\int (S_n(\bold{t}))^2 W(\bold{t})\, d\bold{t},
\end{eqnarray*}
where 
\begin{eqnarray*}
S_n(\bold{t}) &=& \sqrt{n}\sum_{j=1}^n\wj\cos(t_0\hat\eps_j)\Bigg[\cos\Big(\sum_{\nu=1}^k t_\nu X_{j-\nu}\Big)+\sin\Big(\sum_{\nu=1}^k t_\nu X_{j-\nu}\Big)\\
&&\qquad\qquad\qquad\qquad-\sum_{\ell=1}^n\wl\Big(\cos\Big(\sum_{\nu=1}^k t_\nu X_{\ell-\nu}\Big)+\sin\Big(\sum_{\nu=1}^k t_\nu X_{\ell-\nu}\Big)\Big)\Bigg]\\
&&{}+\sqrt{n}\sum_{j=1}^n\wj\sin(t_0\hat\eps_j)\Bigg[\cos\Big(\sum_{\nu=1}^k t_\nu X_{j-\nu}\Big)-\sin\Big(\sum_{\nu=1}^k t_\nu X_{j-\nu}\Big)\\
&&\qquad\qquad\qquad\qquad-\sum_{\ell=1}^n\wl\Big(\cos\Big(\sum_{\nu=1}^k t_\nu X_{\ell-\nu}\Big)-\sin\Big(\sum_{\nu=1}^k t_\nu X_{\ell-\nu}\Big)\Big)\Bigg].
\end{eqnarray*}
For simplicity for the moment we consider only
\begin{eqnarray*}
S_n^{(1)}(\bold{t}) &=& \sqrt{n}\sum_{j=1}^n\wj\cos(t_0\hat\eps_j)\Bigg[\cos\Big(\sum_{\nu=1}^k t_\nu X_{j-\nu}\Big)-\sum_{\ell=1}^n\wl\cos\Big(\sum_{\nu=1}^k t_\nu X_{\ell-\nu}\Big)\Bigg].
\end{eqnarray*}
By a second order Taylor expansion for
\begin{eqnarray*}
\cos\left(t_0\he_j\right)
&=&\cos \left(t_0\left(\e_j+\e_j\frac{\sigma-\hat{\sigma}}{\hat{\sigma}}(X_{j-1})+\frac{m-\hat{m}}{\hat{\sigma}}(X_{j-1})\right)\right)
%&=& \cos\left(t_0\e_j\right)
%-\sin\left(t_0\e_j\right)t_0\left(\e_j\frac{\sigma-\hat{\sigma}}{\hat{\sigma}}(X_{j-1})+\frac{m-\hat{m}}{\hat{\sigma}}(X_{j-1})\right)\\
%&&-\frac 12\cos\left(t_0\xi_j\right)t_0^2\left(\e_j\farc{\sigma-\hat{\sigma}}{\hat{\sigma}}(X_{j-1})+\frac{m-\hat{m}}{\hat{\sigma}}(X_{j-1})\right)^2
\end{eqnarray*}
and introducing the notations
\begin{eqnarray}
\hat{\kappa}_n&=&\frac 1n\sum_{i=1}^nw_n(X_{i-1})\label{kappan}\\
Y_j(\bold{t})&=& \cos\Big(\sum_{\nu=1}^k t_\nu X_{j-\nu}\Big), \quad j=1,\dots,n,\nonumber
\end{eqnarray}
we obtain the expansion $S_n^{(1)}= S_n^{(1,1)}+ S_n^{(1,2)}-\frac 12S_n^{(1,3)}$, where
\begin{eqnarray*}
  S_n^{(1,1)}(\bold{t})&=& \frac{1}{\hat\kappa_n}\frac{1}{\sqrt{n}}\sum_{j=1}^nw_n(X_{j-1})\cos(t_0\eps_j)\Big(Y_j(\bold{t})-\frac{1}{\hat\kappa_n}\frac{1}{n}\sum_{\ell=1}^n w_n(X_{\ell-1})Y_\ell(\bold{t})\Big)\\
  &=& \frac{1}{\hat\kappa_n}\frac{1}{\sqrt{n}}\sum_{j=1}^nw_n(X_{j-1})\Big(\cos(t_0\eps_j)-E[\cos(t_0\eps_j)]\Big)\Big(Y_j(\bold{t})-\frac{1}{\hat\kappa_n}\frac{1}{n}\sum_{\ell=1}^n w_n(X_{\ell-1})Y_\ell(\bold{t})\Big)\\
  S_n^{(1,2)}(\bold{t})&=& \frac{1}{\hat\kappa_n}\frac{1}{\sqrt{n}}\sum_{j=1}^nw_n(X_{j-1})\sin(t_0\e_j)t_0
  \left(\frac{\hat m-{m}}{\hat{\sigma}}(X_{j-1})+\e_j\frac{\hat\sigma-{\sigma}}{\hat{\sigma}}
(X_{j-1})\right)\\
&&{}\times
\Big[Y_j(\bold{t})-\frac{1}{\hat\kappa_n}\frac{1}{n}\sum_{\ell=1}^n w_n(X_{\ell-1})Y_\ell(\bold{t})\Big]\\
  S_n^{(1,3)}(\bold{t})&=& \frac{1}{\hat\kappa_n}\frac{1}{\sqrt{n}}\sum_{j=1}^nw_n(X_{j-1})\cos(t_0\xi_j) t_0^2
  \left(\frac{\hat m-{m}}{\hat{\sigma}}(X_{j-1})+\e_j\frac{\hat\sigma-{\sigma}}{\hat{\sigma}}
(X_{j-1})\right)^2\\
  &&{}\times \Big(Y_j(\bold{t})-\frac{1}{\hat\kappa_n}\frac{1}{n}\sum_{\ell=1}^n w_n(X_{\ell-1})Y_\ell(\bold{t})\Big)
\end{eqnarray*}
(with  $\xi_j$ between $\e_j$ and $\he_j$, $j=1,\dots,n$).
The last term is negligible because
$$\int (S_n^{(1,3)}(\bold{t}))^2W(\bold{t})\,d\bold{t}\leq \int t_0^4W(\bold{t})\,d\bold{t}\Big(\frac{1}{n}\sum_{i=1}^n \eps_i^2+1\Big)^2nO_P\Big(\Big(\frac{a_n^*}{\Delta_n}\Big)^4\Big)=o_P(1)$$
by assumptions \ref{A2} and \ref{A4}, Proposition \ref{rates} and (\ref{kappa}).
Lemmata  \ref{sn11-neu}, \ref{sn12-neu} and \ref{sn12} give further expansions of $S_n^{(1,1)}$ and $S_n^{(1,2)}$. With this we obtain altogether  that $\int (S_n^{(1)}(\bold{t})-\tilde S_n^{(1)}(\bold{t}))^2 W(\bold{t})\, d\bold{t}=o_P(1)$, where
\begin{eqnarray*}
\tilde S_n^{(1)}(\bold{t}) &=& \frac{1}{\sqrt{n}}\sum_{j=1}^n\Bigg[w_n(X_{j-1})\Big(\cos(t_0\eps_j)-E[\cos(t_0\eps_j)]\Big)(Y_j(\bold{t})-E[Y_j(\bold{t})])\\
&& {}\qquad+t_0\Big(\frac{1}{2}(\eps_j^2-1)E[\sin(t_0\eps_1)\eps_1]+\eps_j E[\sin(t_0\eps_1)]\Big)\Big(E[Y_j(\bold{t})|X_{j-1}]-E[Y_j(\bold{t})]\Big)\Bigg].
\end{eqnarray*}
Analogously it follows that $S_n=\tilde S_n+R_n$, where $\int R_n^2(\bold{t}) W(\bold{t})\, d\bold{t}=o_P(1)$ and
\begin{eqnarray*}
&&\tilde S_n(\bold{t})\\
 &=& \frac{1}{\sqrt{n}}\sum_{j=1}^nw_n(X_{j-1})\Bigg[\Big(\cos(t_0\eps_j)-E\big[\cos(t_0\eps_j)\big]\Big)\Big(Y_j(\bold{t})+Z_j(\bold{t})-E[Y_j(\bold{t})+Z_j(\bold{t})]\Big)\\
&& {} +\Big(\sin(t_0\eps_j)-E[\sin(t_0\eps_j)]\Big)\Big(Y_j(\bold{t})-Z_j(\bold{t})-E[Y_j(\bold{t})-Z_j(\bold{t})]\Big)
\\
&& {}+t_0\Big(\eps_j E[\sin(t_0\eps_1)]+\frac{1}{2}(\eps_j^2-1)E[\sin(t_0\eps_1)\eps_1]\Big)\Big(E[Y_j(\bold{t})+Z_j(\bold{t})|X_{j-1}]-E[Y_j(\bold{t})+Z_j(\bold{t})]\Big)\\
&& {} -t_0\Big(\eps_j E[\cos(t_0\eps_1)]+\frac{1}{2}(\eps_j^2-1)E[\cos(t_0\eps_1)\eps_1]\Big)\Big(E[Y_j(\bold{t})-Z_j(\bold{t})|X_{j-1}]-E[Y_j(\bold{t})-Z_j(\bold{t})]\Big)\Bigg],
\end{eqnarray*}
with
$$Z_j(\bold{t})=\sin\Big(\sum_{\nu=1}^k t_\nu X_{j-\nu}\Big)\quad j=1,\dots,n.$$
  To finish the proof  of Theorem \ref{theo1} we apply  Theorem 22
(pages 380, 381) in Ibragimov and Chasminskij (1981). In order to verify the assumptions  
it suffices to show:
 \begin{itemize}
 \item (i) $\tilde S_n(\bold{t})$ has asymptotically normal distribution with zero mean and finite variance;

     \item (ii) for any compact set $F$ in $\mathbb{R}^{k+1}$, 
    $$ \sup_n E\int_F \tilde S^2_n(\bold{t}) W(\bold{t}) dt<\infty;
     $$

     \item (iii)
     $$
     E| \tilde S^2_n(\bold{t}_1)-\tilde S^2_n(\bold{t}_2)|\leq D |\bold{t}_1-\bold{t}_2|^{\gamma}\quad \forall \bold{t}_1, \bold{t}_2
     $$
      for some $\gamma>0$ and some $D>0$;
      
      \item (iv) for all $\eta>0$ there exists some compact set $F_\eta$ in $\mathbb{R}^{k+1}$ with $$E\int_{\mathbb{R}^{k+1}\setminus F_\eta} \tilde S^2_n(\bold{t}) W(\bold{t}) dt<\eta \; \forall n, \; E\int_{\mathbb{R}^{k+1}\setminus F_\eta} \tilde S^2(\bold{t}) W(\bold{t}) dt<\eta.$$
 \end{itemize}
  Since $\tilde S_n(\bold{t})$ is the sums of martingale differences for each $t$ and the  the central limit theorem for martingale differences can be applied which further  implies (i). Direct calculations gives (ii).   Concerning (iii) we have
 \begin{align*}
E &|\tilde S^2_n(\bold{t}_1)-\tilde S^2_n(\bold{t}_2)|\leq E\Big[|\tilde S_n(\bold{t}_1)-\tilde S_n(\bold{t}_2)|\times \Big(|\tilde S^2_n(\bold{t}_1)|+ |\tilde S_n(\bold{t}_2)|\Big)\Big]\\
&\leq \Big(  E\Big|\tilde S_n(\bold{t}_1)-\tilde S_n(\bold{t}_2)\Big|^2\times
E\Big(|\tilde S_n(\bold{t}_1)|+ |\tilde S_n(\bold{t}_2)|\Big)^2\Big)^{1/2}
\end{align*}
 and since
  $$
  E(\tilde S_n(\bold{t})^2)\leq D,\quad \forall \bold{t}
  $$
  it suffices to study
  $$
  E\Big|\tilde S_n(\bold{t}_1)-\tilde S_n(\bold{t}_2)\Big|^2.
  $$
  We show here the needed inequality only for one of  the  terms  in $\tilde S_n(\bold{t}_1)-\tilde S_n(\bold{t}_2)$ all others are treated in the same way. Particularly,
   \begin{align*}
  E&\Big(\frac{1}{\sqrt n}\sum_{j=1}^n w_n(X_{j-1} ) \Big((\cos(t_{01} \varepsilon_j)-
  E\cos(t_{01} \varepsilon_j))(Y_j(\bold{t}_1)-EY_j(\bold{t}_1))\\
  &\qquad\qquad- (\cos(t_{02} \varepsilon_j)-
  E\cos(t_{02} \varepsilon_j))(Y_j(\bold{t}_2)-EY_j(\bold{t}_2))\Big)\Big)^2\\
  &  = E\Big(   w_n(X_{j-1} ) \Big((\cos(t_{01} \varepsilon_j)-
  E\cos(t_{01} \varepsilon_j))(Y_j(\bold{t}_1)-EY_j(\bold{t}_1))\\
&\qquad\qquad- (\cos(t_{02} \varepsilon_j)-
  E\cos(t_{02} \varepsilon_j))(Y_j(\bold{t}_2)-EY_j(\bold{t}_2))\Big)\Big)^2
  \\
  & \leq D || \bold{t}_1-\bold{t}_2||^2
  \end{align*}
   where we used smoothness of cosine and moment assumptions. 
    Proceeding similarly with other terms and putting all together we conclude
    $$
    E |\tilde S^2_n(\bold{t}_1)-\tilde S^2_n(\bold{t}_2)|\leq D || \bold{t}_1-\bold{t}_2||
$$
 This implies the  item (iii).
 Item (iv) follows straightforwardly by our moment assumptions and integrability of $W$.

  \noindent Combining all the above arguments we can infer that the assertion of Theorem   \ref{theo1} holds true;  see Lemma 7.1 and proof of Theorem 4.1 (a) in Hlßvka et al.\ (2014) for a similar argumentation. 
   \hfill $\Box$

\medskip

{\bf Proof of Theorem \ref{theo2}. }

We use the same decomposition of $T_n=\int (S_n(\bold{t}))^2 W(\bold{t})\, d\bold{t}$ as in the proof of Theorem \ref{theo1}. Please note that Lemma \ref{rates} remains true under the assumptions of Theorem \ref{theo2}. A careful inspection of the proof of Theorem \ref{theo1} shows that applying this Lemma one obtains
 $S_n=\tilde S_n+R_n$, where $\int R_n^2(\bold{t}) W(\bold{t})\, d\bold{t}=o_P(n)$ and
\begin{eqnarray*}
\frac{\tilde S_n(\bold{t})}{\sqrt{n}} &=& \frac{1}{n}\sum_{j=1}^nw_n(X_{j-1})\Bigg[\Big(\cos(t_0\eps_j)-E\big[\cos(t_0\eps_j)\big]\Big)\Big(Y_j(\bold{t})+Z_j(\bold{t})-E[Y_j(\bold{t})+Z_j(\bold{t})]\Big)\\
&& {} +\Big(\sin(t_0\eps_j)-E[\sin(t_0\eps_j)]\Big)\Big(Y_j(\bold{t})-Z_j(\bold{t})-E[Y_j(\bold{t})-Z_j(\bold{t})]\Big)\Bigg].
\end{eqnarray*}
 The proof is finished as the end of the proof of Theorem \ref{theo1} applying Theorem 22 (pages 380, 381) in Ibragimov and
Chasminskij (1981). To this end, condition (i) is replaced by convergence in probability of $\tilde S_n(\bold{t})/n^{1/2}$ to \beq
\bar{S}(\bold{t})
&=& E\Bigg[\Big(\cos(t_0\eps_j)-E\big[\cos(t_0\eps_j)\big]\Big)\Big(Y_j(\bold{t})+Z_j(\bold{t})-E[Y_j(\bold{t})+Z_j(\bold{t})]\Big)\\
&& {} +\Big(\sin(t_0\eps_j)-E[\sin(t_0\eps_j)]\Big)\Big(Y_j(\bold{t})-Z_j(\bold{t})-E[Y_j(\bold{t})-Z_j(\bold{t})]\Big)\Bigg]
\eeq
for all $\bold{t}$, whereas in conditions (ii)--(iv) $\tilde S_n$ is replaced by $\tilde S_n/n^{1/2}$. 
Thus we obtain convergence of $T_n/n$ to $\int (\bar{S}(\bold{t}))^2W(\bold{t})\,d\bold{t}$ in probability. 
Note further that by the addition theorems for trigonometric functions and symmetry properties of cosine and sine it holds $\tilde T=\int (\bar{S}(\bold{t}))^2W(\bold{t})d\bold{t}$.
This completes the proof.
\hfill $\Box$

\bigskip

%%%
\noindent{\bf Proof of Lemma \ref{fourier}. }\\
Using assumption \ref{A5}, it follows that $\int_{\mathds R^{k+1}} \vert t_0^4 W(t_0,\ldots,t_k)\vert d(t_0,\ldots,t_k)<\infty$ and since \linebreak $W(t_0,\ldots,t_k)=V_0(t_0)\prod_{i=1}^kV_i(t_i)$ by assumption, one obtains
\begin{align}
\infty>
\int_{\mathds R^{k+1}} \vert t_0^4 W(t_0,\ldots,t_k)\vert d(t_0,\ldots,t_k)
=&
\int_{\mathds R^{k+1}} t_0^4 V_0(t_0)\prod_{i=1}^kV_i(t_i)d(t_0,\ldots,t_k)\notag\\
=&
\int_{\mathds R} t_0^4 V_0(t_0)dt_0
\prod_{i=1}^k\int_{\mathds R} V_i(t_i)dt_i\notag
\end{align}
which gives that $V_i\in L^1(\mathds R)$ for any $i=0,\ldots,k$. Hence, the Fourier transformation of any $V_i$, say $\mathcal F[V_i]$, exists.
The representation of the test statistic is now straightforwardly computed by using the definition of the Fourier transformation and  of the empirical characteristic functions besides the multiplicative structure of $W$.
Since the computation is tedious but without further insights, this part of the proof is omitted here.
\hfill $\Box$

\section{Auxiliary results}
\def\theequation{B.\arabic{equation}}
\setcounter{equation}{0}

First note that for $\hat\kappa_n$ defined in (\ref{kappan}) one obtains directly that $E[(\hat\kappa_n-1)^2]$ can be bounded by $1-F_{X_1}(\frac{a_n}{2})-F_{X_1}(-\frac{a_n}{2})=o(1)$ and thus we have
\begin{equation}\label{kappa}
\hat\kappa_n=1+o_P(1).
\end{equation}

\begin{prop}\label{rates} Let $(X_j)_{j\in \mathbb{Z}}$ be  a strictly stationary time series with marginal density $f_X$. Define $m(x)=E[X_{j+1}\mid X_j=x]$ and $\sigma^2(x)=\Var(X_{j+1}\mid X_j=x)$ and assume \ref{A1}, \ref{Am}, \ref{A3}, \ref{A4}. Let $\Delta_n=\inf_{|x|\leq a_n} f_X(x)$, $a_n^*=((\log n)/(nc_n))^{1/2} +(\log n)^Dc_n^2$ and  $b_n^*= ((c_n/n)^{1/2}+c_n^2)(\log n)^{D}$. Here, $D>0$ is some multiple of $r$ from assumption \ref{Am} and may differ from line to line. We then have
\begin{itemize}
\item[(i)] 
\begin{align*}
\sup_{|x|\leq a_n}|&\hat f_X(x)-f_X(x)| =O_P(a_n^*)\\
\sup_{|x|\leq a_n}|&\hat m(x)-m(x)| =O_P(\frac{a_n^*}{\Delta_n})\\
\sup_{|x|\leq a_n}|&\hat \sigma(x)-\sigma(x)| =O_P(\frac{a_n^*}{\Delta_n})
\end{align*}
\item[(ii)]
\begin{align*}
\sup_{|x|\leq a_n} \Big|& \frac{1}{nc_n}\sum_{j=1}^nK\left(\frac{X_{j-1}-x}{c_n}\right)(m(X_{j-1})-m(x))\Big|=
O_P(b_n^*)\\
\sup_{|x|\leq a_n} \Big|& \frac{1}{nc_n}\sum_{j=1}^nK\left(\frac{X_{j-1}-x}{c_n}\right)(\sigma^2 (X_{j-1})- \sigma^2(x))|=
O_P(b_n^*)\\
\sup_{|x|\leq a_n} \Big|& \frac{1}{nc_n}\sum_{j=1}^nK\left(\frac{X_{j-1}-x}{c_n}\right)(m^2 (X_{j-1})- m^2(x))|=
O_P(b_n^*)
\end{align*}
\item[(iii)]
\begin{align*}
\sup_{|x|\leq a_n} \Big|& \frac{1}{nc_n}\sum_{j=1}^nK\left(\frac{X_{j-1}-x}{c_n}\right)\sigma(X_{j-1})\e_j(m (X_{j-1})- m(x))|=
O_P(b_n^*).
\end{align*}
\end{itemize}
\end{prop}

{\bf Proof.} The first two results of (i) are stated in Theorems 6 and 8 by Hansen (2008)  without the $(\log n)^D$ factor of the $c_n^2$ term. 
In comparison to Hansen (2008) we use a different bounding for the expectation terms since we do not assume second derivatives. E.\,g.\ we obtain, making use of the mean value theorem, the properties of the kernel function and our assumption \ref{Am},
\begin{eqnarray*}
\sup_{|x|\leq a_n}|E[\hat f_X(x)-f_X(x)]|&=&\sup_{|x|\leq c_n}|\int K(u)(f_X(x-c_nu)-f_X(x))\,du|\\
&\leq& \sup_{|x|\leq a_n} c_n\int K(u)|u|\sup_{\xi \text{ between } \atop x \text{ and }x-c_nu}|f'(\xi)-f'(x)|\,du\\
&=& O(c_n^2 (\log n)^r).
\end{eqnarray*}
 The result on $\hat\sigma$ follows similarly to the derivations by Hansen (2008) by noting that $\hat\sigma^2(x)=\hat s(x)-\hat m^2(x)$, where $\hat s$ is the Nadaraya-Watson estimator for $s(x)=E[X_j^2\mid X_j=x]$ based on the observation pairs $(X_{j-1},X_j^2)$, $j=1,\dots,n$.

Towards  the results in (ii)  we treat only the first one since the others follow analogously.  Note that by the mean value theorem 
 $$
 m(X_{i-1})-m(x)=(X_{i-1}-x) m'(x)+(X_{i-1}-x) (m'(\xi_{X_{i-1},x})-m'(x))
 $$
 for some $\xi_{X_{i-1},x}$ between $[\min(X_{i-1},x),\max(X_{i-1}, x)]$, where the absolute value of the second summand can be bounded by $(X_{i-1}-x)^2(\log n)^r$ due to assumption \ref{Am}. It thus suffices to show
 \begin{align*}
\sup_{|x|\leq a_n}\Big| \frac{1}{nc_n}\sum_{i=1}^nK\left(\frac{X_{i-1}-x}{c_n}\right)(X_{i-1}- x)^2 \Big|  &= O_P(c_n^2)
\\
 \sup_{|x|\leq a_n}\Big| \frac{1}{nc_n}\sum_{i=1}^nK\left(\frac{X_{i-1}-x}{c_n}\right)(X_{i-1}-x) m'(x)\Big|  &= O_P(b_n^*).
 \end{align*}
The first relation is straightforward by assumption \ref{Am}  and applying  Theorem 2 in Hansen (2008) with $Y_i=1$ and the kernel $u\mapsto K(u)u^2$. For the latter one we receive with the same theorem applied with $Y_i=1$ and kernel $u\mapsto K(u)u$
  \begin{align*}
  \sup_{|x|\leq a_n}&\Big| \frac{1}{nc_n}\sum_{i=1}^n\Big(K\left(\frac{x-X_{i-1}}{c_n}\right)
             (x-X_{i-1})- E\Big[ K\left(\frac{x-X_{i-1}}{c_n}\right)
             (x-X_{i-1})\Big]\Big)\Big|
             \\
             &
             =O_P\Big(\Big(\frac{\log n}{nc_n}\Big)^{1/2}c_n\Big).
\end{align*}
Further by direct calculation
   \begin{align*}
    \frac 1{nc_n}\sum_{i=1}^n &E\Big[ K\left(\frac{x-X_{i-1}}{c_n}\right)
             (x-X_{i-1})\Big]= \frac{1}{c_n}   \int   K\left(\frac{x-y}{c_n}\right) (x-y) f_{X_{i-1}}(y) dy\\
      &=c_n\int K(z) z f_{X_{i-1}}(x-zc_n)) dz=O( c_n^2 \int z^2 K(z)dz \sup_{x\in\tilde I_n}| f'_{X}(x)|)=O( c_n^2 (\log n)^{r})
    \end{align*}
         where we utilize  assumptions \ref{Am} and \ref{A3}.

 The result (iii) can be proved in the same way as the results in (ii). Just set $Y_i=|\e_i|$ for the first and $Y_i=\sigma(X_{i-1})\e_i$ for the second relation (when applying Theorem 2 by Hansen, 2008) and note that $$E\Big[ K\left(\frac{x-X_{i-1}}{c_n}\right)(x-X_{i-1})\sigma(X_{i-1})\e_i\Big]=0.$$
       \hfill $\Box$

\begin{lemma}\label{sn11-neu}  Under the assumptions of Theorem \ref{theo1} we have $$\int (S_n^{(1,1)}(\bold{t})-\tilde S_n^{(1,1)}(\bold{t}))^2 W(\bold{t})\, d\bold{t}=o_P(1),$$ where
$$\tilde S_n^{(1,1)}(\bold{t})=\frac 1{\sqrt n}\sum_{j=1}^nw_n(X_{j-1})\big(\cos(t_0\eps_j)-E[\cos(t_0\e_j)]\big)(Y_j(\bold{t})-E[Y_j(\bold{t})]).$$
\end{lemma}

{\bf Proof.} Due to (\ref{kappa}) we have $$ S_n^{(1,1)}(\bold{t})=(1+o_P(1))(\tilde S_n^{(1,1)}(\bold{t})-J_n(\bold{t})J_n^{(1)}(\bold{t})-J_n(\bold{t})J_n^{(2)}(\bold{t})),$$ where
\begin{eqnarray*}
J_n(\bold{t}) &=& \frac{1}{\sqrt{n}}\sum_{j=1}^nw_n(X_{j-1})\Big(\cos(t_0\eps_j)-E[\cos(t_0\eps_j)]\Big)\\
 J_n^{(1)}(\bold{t})&=& \frac{1}{n}\sum_{\ell=1}^n \Big(Y_\ell(\bold{t})-E[Y_\ell(\bold{t})]\Big)\\
 J_n^{(2)}(\bold{t})&=& \frac{1}{n}\sum_{\ell=1}^n  Y_\ell(\bold{t}) (w_n(X_{\ell-1})-\hat\kappa_n)\frac{1}{\hat\kappa_n}.
\end{eqnarray*}
Note that from  assumption \ref{A1} it follows that $\beta>2$ and thus $\sum_{i=0}^\infty (i+1)\alpha(i)<\infty$. From this, centeredness of the summands (under the null) and the boundedness of cosine analogously to the proof of Theorem 2 by Yokoyama (1980) one obtains
\begin{eqnarray}
&&E[(J_n(\bold{t}))^4] \leq D\label{yoko-0}\\
&&E[(J_n^{(1)}(\bold{t}))^4] \leq  \frac{1}{n^2}D.\label{yoko}
\end{eqnarray}
The constant $D$ can be chosen independent of $\bold{t}$ due to the boundedness of the cosine function. Thus from the Cauchy Schwarz inequality we obtain directly
\begin{eqnarray*}
E\Big[\int ( J_n(\bold{t})J_n^{(1)}(\bold{t}))^2W(\bold{t})\,d\bold{t}\Big] = O(\frac{1}{n}).
\end{eqnarray*}
Now note that
\begin{eqnarray*}
 J_n^{(2)}(\bold{t}) %&=& \frac{1}{n}\sum_{\ell=1}^n  (Y_\ell(t)-E[Y_\ell(t)]) (w_n(X_{\ell-1})-\hat\kappa_n)\frac{1}{\hat\kappa_n}\\
%&=& (1+o_P(1))\Bigg(\frac{1}{n}\sum_{\ell=1}^n  (Y_\ell(t)-E[Y_\ell(t)]) (w_n(X_{\ell-1})-E[w_n(X_{\ell-1})])\\
%&&{}-\frac{1}{n}\sum_{\ell=1}^n  (Y_\ell(t)-E[Y_\ell(t)])\frac{1}{n}\sum_{j=1}^n (w_n(X_{j-1})-E[w_n(X_{j-1})])\Bigg)
&=& (1+o_P(1))\Bigg(\frac{1}{n}\sum_{\ell=1}^n  Y_\ell(\bold{t}) (w_n(X_{\ell-1})-E[w_n(X_{\ell-1})])\\
&&{}-\frac{1}{n}\sum_{\ell=1}^n  Y_\ell(\bold{t})\frac{1}{n}\sum_{j=1}^n (w_n(X_{j-1})-E[w_n(X_{j-1})])\Bigg)
\end{eqnarray*}
and thus by boundedness of $Y_\ell$ we have, uniformly with respect to $\bold{t}$,
\begin{eqnarray}\label{j2}
| J_n^{(2)}(\bold{t})|
&=& O_P(1)\frac{1}{n}\sum_{j=1}^n |w_n(X_{j-1})-E[w_n(X_{j-1})]|=o_P(1)
\end{eqnarray}
by a consideration of the expectation of the sum  due to the properties of the weight function. We obtain
\begin{eqnarray*}
\int ( J_n(\bold{t})J_n^{(2)}(\bold{t}))^2W(\bold{t})\,d\bold{t} &=& o_P(1)\int ( J_n(\bold{t}))^2W(\bold{t})\,d\bold{t} \;=\; o_P(1)
\end{eqnarray*}
by an application of (\ref{yoko-0}).
\hfill $\Box$

\begin{lemma}\label{sn12-neu}  Under the assumptions of Theorem \ref{theo1} we have $$\int (S_n^{(1,2)}(\bold{t})-\tilde S_n^{(1,2)}(\bold{t}))^2 W(\bold{t})\, d\bold{t}=o_P(1),$$ where
$$\tilde S_n^{(1,2)}(\bold{t})=\frac 1{\sqrt n}\sum_{j=1}^nw_n(X_{j-1})\sin(t_0\e_j)t_0\left(\frac{\hat m-{m}}{\hat{\sigma}}(X_{j-1})+\e_j\frac{\hat\sigma-{\sigma}}{\hat{\sigma}}(X_{j-1})\right)(Y_j(\bold{t})-E[Y_j(\bold{t})]).$$
\end{lemma}

{\bf Proof.} Due to (\ref{kappa}) we have $$ S_n^{(1,2)}(\bold{t})=(1+o_P(1))(\tilde S_n^{(1,2)}(\bold{t})-t_0I_n(\bold{t})J_n^{(1)}(\bold{t})-t_0I_n(\bold{t})J_n^{(2)}(\bold{t}))$$ with $J_n^{(1)}$ and $J_n^{(2)}$ as in Lemma \ref{sn11-neu} and
\begin{eqnarray*}
 I_n(\bold{t})&=& \frac 1{\sqrt n}\sum_{j=1}^nw_n(X_{j-1})\sin(t_0\e_j)\left(\frac{\hat m-{m}}{\hat{\sigma}}(X_{j-1})+\e_j\frac{\hat\sigma-{\sigma}}{\hat{\sigma}}(X_{j-1})\right).
\end{eqnarray*}
Now
\begin{eqnarray*}
 |I_n(\bold{t})|&\leq& \frac 1{\sqrt n}\sum_{j=1}^n(|\e_j|+1)O_P(\frac{a_n^*}{\Delta_n}) =o_P(\sqrt n)
\end{eqnarray*}
uniformly with respect to $t$ by assumption \ref{A2} and Proposition \ref{rates} (i). Thus
\begin{eqnarray*}
\int (t_0I_n(\bold{t})J_n^{(1)}(\bold{t}))^2W(\bold{t})\,d\bold{t} &=& o_P( n)\int t_0^2 ( J_n^{(1)}(\bold{t}))^2W(\bold{t})\,d\bold{t} \;=\; o_P(1)
\end{eqnarray*}
by (\ref{yoko}).

Further, by (\ref{j2}) we obtain
\begin{eqnarray*}
\int (t_0I_n(\bold{t})J_n^{(2)}(\bold{t}))^2W(\bold{t})\,d\bold{t} &=& o_P(1)\int t_0^2 ( I_n(\bold{t}))^2W(\bold{t})\,d\bold{t} \;=\; o_P(1),
\end{eqnarray*}
where one yields the last equality as follows.
Similarly to the proof of Lemma \ref{sn12} one can first replace the random denominators $\hat\sigma\hat f_X$ in the definition of $I_n$ by their true counterparts $\sigma f_X$ applying Proposition \ref{rates}. Let $\tilde I_n$ denote the resulting term, then $E[\int t_0^2 ( \tilde I_n(\bold{t}))^2W(\bold{t})\,d\bold{t}]=O(1)$ is shown by straightforward calculations.
\hfill $\Box$

\begin{lemma}\label{sn12} Under the assumptions of Theorem \ref{theo1} we have $$\int \left(\tilde S_n^{(1,2)}(\bold{t})-\tilde S_n^{(1,2,1)}(\bold{t})-\tilde S_n^{(1,2,2)}(\bold{t})\right)^2 W(\bold{t})\,d\bold{t}=o_P(1),$$ where
\begin{eqnarray*}
\tilde S_n^{(1,2,1)}(\bold{t}) &=& \frac{t_0}{\hk n^{3/2}}\sum_{j=1}^n\sum_{i=1}^n \frac{1}{c_n}K\left(\frac{X_{j-1}-X_{i-1}}{c_n}\right)\frac{w_n(X_{j-1})\sigma ( X_{i-1})\eps_i}{\sigma(X_{j-1})f_X(X_{j-1})}\sin(t_0\eps_j)\\
&&\qquad\qquad\qquad\times(Y_j(\bold{t})-E[Y_j(\bold{t})])\\
\tilde S_n^{(1,2,2)}(\bold{t}) &=& \frac{t_0}{\hk n^{3/2}}\sum_{j=1}^n\sum_{i=1}^n \frac{1}{c_n}K\left(\frac{X_{j-1}-X_{i-1}}{c_n}\right)\frac{w_n(X_{j-1})\sigma^2 ( X_{i-1})(\eps_i^2-1)}{2\sigma^2(X_{j-1})f_X(X_{j-1})}\sin(t_0\eps_j)\eps_j\\
&&\qquad\qquad\qquad\times(Y_j(\bold{t})-E[Y_j(\bold{t})])
\end{eqnarray*}

\end{lemma}

{\bf Proof.}
Recall the definition of $\tilde S_n^{(1,2)}$ in Lemma \ref{sn12-neu}  and note that
\begin{eqnarray*}
&&\frac {t_0}{\sqrt n}\sum_{j=1}^n\frac{w_n(X_{j-1})}{\hk}\sin(t_0\eps_j)\frac{ m-\hat m}{\hat{\sigma}}(X_{j-1})(Y_j(\bold{t})-E[Y_j(\bold{t})])\\
&=& \frac {t_0}{\sqrt n}\sum_{j=1}^n\frac{w_n(X_{j-1})}{\hk}\sin(t_0\eps_j)\frac{m-\hat{m}}{\sigma}(X_{j-1})\frac{\hat{f}_X}{f_X}(X_{j-1})(Y_j(\bold{t})-E[Y_j(\bold{t})])
+ R_n^{(1)}(\bold{t})\\
&&{} +R_n^{(2)}(\bold{t})\\
&=& \tilde S_n^{(1,2,1)}(\bold{t})+ R_n^{(1)}(\bold{t}) +R_n^{(2)}(\bold{t})+R_n^{(3)}(\bold{t}),
\end{eqnarray*}
where
\begin{eqnarray*}
 R_n^{(1)}(\bold{t}) &=& \frac {t_0}{\sqrt n}\sum_{j=1}^n\frac{w_n(X_{j-1})}{\hk}\sin(t_0\eps_j)\frac{\hat m-{m}}{\hat{\sigma}}(X_{j-1})\frac{\sigma-\hat{\sigma}}{\sigma}(X_{j-1})(Y_j(\bold{t})-E[Y_j(\bold{t})])
\\
R_n^{(2)}(\bold{t})&=& \frac {t_0}{\sqrt n}\sum_{j=1}^n\frac{w_n(X_{j-1})}{\hk}\sin(t_0\eps_j)\frac{\hat m-{m}}{\sigma}(X_{j-1})\frac{f_X-\hat{f}_X}{f_X}(X_{j-1})(Y_j(\bold{t})-E[Y_j(\bold{t})])
\\
R_n^{(3)}(\bold{t})&=& \frac {t_0}{\hk\sqrt n}\sum_{j=1}^n\frac{w_n(X_{j-1})}{\sigma(X_{j-1})}\sin(t_0\eps_j)\left(\frac{\frac 1{nc_n}\sum_{i=1}^nK\left(\frac{X_{j-1}-X_{i-1}}{c_n}\right)(m(X_{i-1})-m(X_{j-1}))}{f_X(X_{j-1})}\right)\\
&&\qquad\qquad\times(Y_j(\bold{t})-E[Y_j(\bold{t})])).
\end{eqnarray*}
By Proposition \ref{rates} (i) one directly obtains that $\int (R_n^{(j)}(\bold{t}))^2W(\bold{t})\,d\bold{t}$ for $j=1,2$ is of rate $O_P(n(a_n^*/\Delta_n)^4)=o_P(1)$.

Concerning  $\int (R_n^{(3)}(\bold{t}))^2W(\bold{t})\,d\bold{t}$  notice that
$$
|R_n^{(3)}(\bold{t})|\leq  D \frac {t_0}{\hk\sqrt n}\sum_{j=1}^n\frac{w_n(X_{j-1})}{\sigma(X_{j-1})f_X(X_{j-1})}\sup_{|x|\leq a_n}
\Big|\frac{ 1}{nc_n}\sum_{i=1}^nK\left(\frac{x-X_{i-1}}{c_n}\right)(m(X_{i-1})-m(x))\Big|
$$
uniformly in $\bold{t}$ which  together with assertion \ref{rates} (ii)  implies the rate  $\int (R_n^{(3)}(\bold{t}))^2W(\bold{t})\,d\bold{t}=O_P(n(\log n)^{4r}(b_n^*)^2)=o_P(1)$, where the latter equality follows from assumption \ref{Aneu}.

Concerning the second term in the definition of $\tilde S_n^{(1,2)}$ in Lemma \ref{sn12-neu} note that due to $\hat\sigma-\sigma=(\hat\sigma^2-\sigma^2)/(\hat\sigma+\sigma)$, analogous to before one shows that
\begin{eqnarray*}
&&\frac {t_0}{\sqrt n}\sum_{j=1}^n\frac{w_n(X_{j-1})}{\hk}\sin(t_0\eps_j)\e_j\frac{\hat\sigma-{\sigma}}{\hat{\sigma}}(X_{j-1})(Y_j(\bold{t})-E[Y_j(\bold{t})])\\
&=& \frac {t_0}{\sqrt n}\sum_{j=1}^n\frac{w_n(X_{j-1})}{\hk}\sin(t_0\eps_j)\e_j\frac{(\hat\sigma^2-{\sigma}^2)\hat f_X}{2\sigma^2 f_X}(X_{j-1})(Y_j(\bold{t})-E[Y_j(\bold{t})])
+R_n^{(4)}(\bold{t}),
\end{eqnarray*}
where $\int (R_n^{(4)}(\bold{t}))^2W(\bold{t})\,d\bold{t}=o_P(1)$ now follows from Proposition \ref{rates} (i). 

To treat the remaining term further we first insert the definition of $\hat\sigma^2$ and then use the fact that $m^2-\hat m^2=2m(m-\hat m)-(m-\hat m)^2$   and insert the definition of $\hat m$.
With this one obtains 
\begin{eqnarray*}
&&\frac {t_0}{\sqrt n}\sum_{j=1}^n\frac{w_n(X_{j-1})}{\hk}\sin(t_0\eps_j)\e_j\frac{(\hat\sigma^2-{\sigma}^2)\hat f_X}{2\sigma^2 f_X}(X_{j-1})(Y_j(\bold{t})-E[Y_j(\bold{t})])\\
&=& \tilde S_n^{(1,2,2)}(\bold{t})+R_n^{(5)}(\bold{t})+R_n^{(6)}(\bold{t})+R_n^{(7)}(\bold{t})+R_n^{(8)}(\bold{t})-R_n^{(9)}(\bold{t}),
\end{eqnarray*}
where
\begin{eqnarray*}
 R_n^{(5)}(\bold{t}) &=& \frac{t_0}{\hk n^{3/2}}\sum_{j=1}^n\sum_{i=1}^n \frac{1}{c_n}K\left(\frac{X_{j-1}-X_{i-1}}{c_n}\right)\frac{w_n(X_{j-1})}{2\sigma^2(X_{j-1})f_X(X_{j-1})}\sin(t_0\eps_j)\eps_j\\
&&\qquad\qquad\times(\sigma^2 ( X_{i-1})-\sigma^2 ( X_{j-1}))(Y_j(\bold{t})-E[Y_j(\bold{t})])
\\
R_n^{(6)}(\bold{t}) &=& \frac{t_0}{\hk n^{3/2}}\sum_{j=1}^n\sum_{i=1}^n \frac{1}{c_n}K\left(\frac{X_{j-1}-X_{i-1}}{c_n}\right)\frac{w_n(X_{j-1})}{2\sigma^2(X_{j-1})f_X(X_{j-1})}\sin(t_0\eps_j)\eps_j\\
&&\qquad\qquad\times(m^2 ( X_{i-1})-m^2 ( X_{j-1}))(Y_j(\bold{t})-E[Y_j(\bold{t})])
\\
R_n^{(7)}(\bold{t}) &=& \frac{t_0}{\hk n^{3/2}}\sum_{j=1}^n\sum_{i=1}^n \frac{1}{c_n}K\left(\frac{X_{j-1}-X_{i-1}}{c_n}\right)\frac{w_n(X_{j-1})}{\sigma^2(X_{j-1})f_X(X_{j-1})}\sin(t_0\eps_j)\eps_j\\
&&\qquad\qquad\times m(X_{j-1})(m ( X_{j-1})-m ( X_{i-1}))(Y_j(\bold{t})-E[Y_j(\bold{t})])
\\
R_n^{(8)}(\bold{t}) &=& \frac{t_0}{\hk n^{3/2}}\sum_{j=1}^n\sum_{i=1}^n \frac{1}{c_n}K\left(\frac{X_{j-1}-X_{i-1}}{c_n}\right)\frac{w_n(X_{j-1})}{\sigma^2(X_{j-1})f_X(X_{j-1})}\sin(t_0\eps_j)\eps_j\\
&&\qquad\qquad\times\sigma(X_{i-1})\e_i(m(X_{i-1})-m(X_{j-1}))(Y_j(\bold{t})-E[Y_j(\bold{t})])
\\
R_n^{(9)}(\bold{t})&=& \frac{t_0}{\hk n^{3/2}}\sum_{j=1}^n\sum_{i=1}^n \frac{1}{c_n}K\left(\frac{X_{j-1}-X_{i-1}}{c_n}\right)\frac{w_n(X_{j-1})}{2\sigma^2(X_{j-1})f_X(X_{j-1})}\sin(t_0\eps_j)\eps_j\\
&&\qquad\qquad\times (m(X_{j-1})-\hat m(X_{j-1}))^2(Y_j(\bold{t})-E[Y_j(\bold{t})])
\end{eqnarray*}
and one can show $\int (R_n^{(j)}(\bold{t}))^2W(\bold{t})\,d\bold{t}=o_P(1)$ completely analogous to the treatment of $R_n^{(3)}$ for $j=5,6,7,8$ and $R_n^{(1)}$ for $j=9$.
\hfill $\Box$

\bigskip

\begin{lemma}\label{sn122} Under the assumptions of Theorem \ref{theo1} we have 
\begin{eqnarray*}
\int (\tilde S_n^{(1,2,1)}(\bold{t})-\bar S_n^{(1,2,1)}(\bold{t}))^2 W(\bold{t})\,d\bold{t}&=&o_P(1)\\
\int (\tilde S_n^{(1,2,2)}(\bold{t})-\bar S_n^{(1,2,2)}(\bold{t}))^2 W(\bold{t})\,d\bold{t}&=&o_P(1),
\end{eqnarray*} 
where
\begin{eqnarray*}
\bar S_n^{(1,2,1)}(\bold{t}) &=& \frac{t_0E[\sin(t_0\eps_1)]}{\sqrt{n}}\sum_{j=1}^n w_n(X_{j-1})\eps_j (E[Y_j(\bold{t})|X_{j-1}]-E[Y_j(\bold{t})])\\
\bar S_n^{(1,2,2)}(\bold{t}) &=& \frac{t_0E[\sin(t_0\eps_1)\eps_1]}{2\sqrt{n}}\sum_{j=1}^n w_n(X_{j-1})(\eps_j^2-1)(E[Y_j(\bold{t})|X_{j-1}]-E[Y_j(\bold{t})]).
\end{eqnarray*}

\end{lemma}

{\bf Proof.}
We only prove the first assertion, the second one can be shown completely analogous.
We have the expansion
\begin{eqnarray*}
\tilde S_n^{(1,2,1)}(\bold{t})-\bar S_n^{(1,2,1)}(\bold{t}) &=& U_n(\bold{t})+t_0E[\sin(t_0\eps_1)]V_n(\bold{t}),
\end{eqnarray*}
where
\begin{eqnarray*}
U_n(\bold{t}) &=& \frac{1}{n^{3/2}}\sum_{j=1}^n\sum_{i=1}^n \varphi(\bold{t},\eps_i,X_{i-1},\zeta_j)
\end{eqnarray*}
with $\zeta_j=(X_{j-1},\dots,X_{j-k})$,
%{\color{red}
% Here I would also  use notation
%$$
%\xi_i=(\eps_i, X_{i-1})
%$$
%Using this notation on the next page we should write
%$$
%g(\xi_{i_1}, \xi_{i_2}, \zeta_{j_1}, \zeta_{j_2})
%$$
%instead of
%$$g(\xi_{i_1}, \xi_{i_2}, \xi_{j_1}, \xi_{j_2})
%$$
\begin{eqnarray*}
&&\varphi(\bold{t},\eps_i,X_{i-1},\zeta_j) \\
&=&
\sigma(X_{i-1})\eps_i\Bigg(
 \frac{1}{c_n}K\Big(\frac{X_{j-1}-X_{i-1}}{c_n}\Big)\frac{w_n(X_{j-1})}{f_X(X_{j-1})\sigma(X_{j-1})}\sin(t_0\eps_j)(Y_j(\bold{t})-E[Y_j(\bold{t})])
\\
&&{}\qquad\qquad -\int \frac{1}{c_n}K\big(\frac{x-X_{i-1}}{c_n}\big)\frac{w_n(x)}{\sigma(x)}E[\sin(t_0\eps_1)](E[Y_j(\bold{t})|X_{j-1}=x]-E[Y_j(\bold{t})])\,dx
\Bigg)
\end{eqnarray*}
and
\begin{eqnarray*}
V_n(\bold{t}) &=& \frac{1}{\sqrt{n}}\sum_{i=1}^n \sigma(X_{i-1})\eps_i\int \frac{1}{c_n}K\big(\frac{x-X_{i-1}}{c_n}\big)
\Big(\frac{ w_n(x)\psi(\bold{t},x)}{\sigma(x)}
-\frac{w_n(X_{i-1})\psi(\bold{t},X_{i-1})}{\sigma(X_{i-1})}\Big)dx
\end{eqnarray*}
with $\psi(\bold{t},x)=E[Y_i(\bold{t})|X_{i-1}=x]-E[Y_i(\bold{t})]$.
Straightforwardly we obtain negligibility of $V_n$ by considering the expectation
\begin{eqnarray*}
&&E\big[\int t_0^2 V_n^2(\bold{t})W(\bold{t})\,d\bold{t}\big] \\
&=& \int t_0^2\int \sigma^2(z)\Big(\int \frac{1}{c_n}K\big(\frac{x-z}{c_n}\big)
\Big(\frac{ w_n(x)\psi(\bold{t},x)}{\sigma(x)}
-\frac{w_n(z)\psi(\bold{t},z)}{\sigma(z)}\Big)dx\Big)^2f_X(z)\,dz\, W(\bold{t})\,d\bold{t}.
\end{eqnarray*}
Note that the inner integral is zero for $z\not\in I_n$. We further separately consider the cases $z\in K_n=[-a_n+c_nC,a_n-c_nC]$ and $z\in I_n\setminus K_n$ to obtain
\begin{eqnarray*}
&&E\big[\int t_0^2 V_n^2(\bold{t})W(\bold{t})\,d\bold{t}\big] \\
&\leq & \int t_0^2\Big[\int \sigma^2(z)\Big(\int \frac{1}{c_n}K\big(\frac{x-z}{c_n}\big)
\Big(\frac{|\psi(\bold{t},x)-\psi(\bold{t},z)|}{\sigma(x)}+\Big|\frac{1}{\sigma(x)}-\frac{1}{\sigma(z)}\Big|\Big)dx\Big)^2\\
&&\qquad\qquad\times f_X(z)I\{z\in K_n\}\,dz\\
&&{}\quad +\int \sigma^2(z)\Big(\int \frac{1}{c_n}K\big(\frac{x-z}{c_n}\big)
\Big(\frac{1}{\sigma(x)}+\frac{1}{\sigma(z)}\Big)\,dx\Big)^2f_X(z)I\{z\in I_n\setminus K_n\}\,dz
\Big] W(\bold{t})\,d\bold{t}\\
&=& O((\log n)^{5r}(c_n+c_n^D)) = o(1)
\end{eqnarray*}
by assumptions \ref{Am} and \ref{Apsi}.

We will now prove $E[\int U_n^2(\bold{t})W(\bold{t})\,d\bold{t}]=o(1)$. To this end note that
\begin{eqnarray}\label{sum_u}
E[\int U_n^2(\bold{t})W(\bold{t})\,d\bold{t}]
&=& \frac{1}{n^{3}}\sum_{j_1=1}^n\sum_{i_1=1}^n\sum_{j_2=1}^n\sum_{i_2=1}^n E\Big[g(\xi_{i_1},\xi_{i_2},\xi_{j_1},\xi_{j_2})\Big],
\end{eqnarray}
where $\xi_i=(\eps_i,\zeta_i)$ and
$$g(\xi_{i_1},\xi_{i_2},\xi_{j_1},\xi_{j_2})=\int \varphi(\bold{t},\eps_{i_1},X_{i_1-1},\zeta_{j_1})
\varphi(\bold{t},\eps_{i_2},X_{i_2-1},\zeta_{j_2})W(\bold{t})\,d\bold{t}$$

We first consider the case where all indices $i_1,j_1,i_2,j_2$ are different. Then the expectation is zero if either $i_1$ or $i_2$ is the largest index because $E[\eps_i]=0$ and $\eps_i$ is independent of $\eps_j,X_{i-1},X_{i-2},\dots$ (for $j\neq i$). All other cases are treated similarly and thus we only discuss the case $i_1 < i_2 <j_1 <j_2$ in detail.
We will apply a version of Lemma 2.1 by Sun and Chiang (1997) for multivariate random variables (see Su and Xiao's  (2008) Lemma D.1) in two separate subcases.
First let $i_2-i_1\geq j_1-i_2$. Denote by the process $\xi_i^*$, $i\in\mathbb{Z}$, an independent copy of $\xi_i$, $i\in\mathbb{Z}$ i.\,e.\ a process with the same distributional properties, but independent of the original data.
Then $E[g(\xi_{i_1}^*,\xi_{i_2},\xi_{j_1},\xi_{j_2})]=0$ and, for $\delta>0$,
\begin{eqnarray*}
&&E\Big[\Big|g(\xi_{i_1}^*,\xi_{i_2},\xi_{j_1},\xi_{j_2})\big|^{1+\delta}\Big]\\
&\leq& k_1 \sup_{x\in [-a_n-Cc_n,a_n+Cc_n]}\sigma^{2+2\delta}(x)\sup_{x\in [-a_n-Cc_n,a_n+Cc_n]}\sigma^{-2-2\delta}(x)E[|\eps_1|^{1+\delta}] \\
&&\times E\Bigg[\int \Bigg( \Big|\frac{1}{c_n}K\Big(\frac{X_{j_1-1}-y}{c_n}\Big)\frac{w_n(X_{j_1-1})}{f_X(X_{j_1-1})}\Big|+\int
\Big|\frac{1}{c_n}K\Big(\frac{x-y}{c_n}\Big)w_n(x)\Big|\,dx\Bigg)^{1+\delta}f_X(y)\,dy\\
&&\quad\times \Bigg( \Big|\frac{1}{c_n}K\Big(\frac{X_{j_2-1}-X_{i_2-1}}{c_n}\Big)\frac{w_n(X_{j_2-1})}{f_X(X_{j_2-1})}\Big|+\int
\Big|\frac{1}{c_n}K\Big(\frac{x-X_{i_2-1}}{c_n}\Big)w_n(x)\Big|\,dx\Bigg)^{1+\delta} |\varepsilon_{i_2}|^{1+\delta}
\Bigg]\\
&\leq& k_2 \sup_{x\in [-a_n-Cc_n,a_n+Cc_n]}\sigma^{2+2\delta}(x)\sup_{x\in [-a_n-Cc_n,a_n+Cc_n]}\sigma^{-2-2\delta}(x)\\
&&{}\times\Big(\sup_{x\in [-a_n-Cc_n,a_n+Cc_n]}(f_X(x))^{-1-\delta}+1\Big)\\
&&{}\times c_n^{-\delta}  E\Bigg[ \Bigg( \Big|\frac{1}{c_n}K\Big(\frac{X_{j_2-1}-X_{i_2-1}}{c_n}\Big)\frac{w_n(X_{j_2-1})}{f_X(X_{j_2-1})}\Big|+1\Bigg)^{1+\delta}|\varepsilon_{i_2}|^{1+\delta}
\Bigg]
\end{eqnarray*}
for some constants $k_1,k_2$.  This is of order $O((\log n)^{\tilde r}c_n^{-2\delta})$ for $\tilde r=5r(1+\delta)$ by assumptions \ref{Am}--\ref{A3}.
%\footnote{we need to specify the assumptions. } %if we do not want to assume joint density (here of $X_{j_2-1},X_{i_2-1}$) we only obtain the rate $$O((\log n)^rc_n^{-1-2\delta})$$
An application of the aforementioned inequality gives
\begin{eqnarray*}
\Bigg|\frac{1}{n^{3}}\sum_{i_1 < i_2 <j_1 <j_2\atop i_2-i_1\geq j_1-i_2} E\Big[g(\xi_{i_1},\xi_{i_2},\xi_{j_1},\xi_{j_2})\Big]\Bigg|
&=& O((\log n)^{{\tilde r}/(1+\delta)})\frac{1}{n^3c_n^{2\delta/(1+\delta)}}\sum_{i_1 < i_2 <j_1 <j_2\atop i_2-i_1\geq j_1-i_2} (\alpha(i_2-i_1))^{\delta/(1+\delta)}\\
&=& O\Big(\frac{(\log n)^{5r}}{nc_n^{2\delta/(1+\delta)}}\Big)\sum_{j=1}^n j(\alpha(j))^{\delta/(1+\delta)} \\
&\leq& o(1)\sum_{j=1}^\infty j^{1-\frac{\beta\delta}{1+\delta}}\;=\; o(1)
\end{eqnarray*}
by assumptions \ref{A1}, \ref{A2} and \ref{Aneu}. Here the mixing coefficient $\alpha$ of $\xi_i$, $i\in \mathbb{Z}$, is the same as the mixing coefficient of $X_i$, $i\in\mathbb{Z}$, see Fan and Yao (2003).
In the subcase $i_2-i_1<j_1-j_2$ we apply the same inequality but by considering $E[g(\xi_{i_1}^*,\xi_{i_2}^*,\xi_{j_1},\xi_{j_2})]=0$ and obtain
\begin{eqnarray*}
\Bigg|\frac{1}{n^{3}}\sum_{i_1 < i_2 <j_1 <j_2\atop i_2-i_1\geq j_1-i_2} E\Big[g(\xi_{i_1},\xi_{i_2},\xi_{j_1},\xi_{j_2})\Big]\Bigg|
&=& O((\log n)^{{\tilde r}/(1+\delta)})\frac{1}{n^3c_n^{2\delta/(1+\delta)}}\sum_{i_1 < i_2 <j_1 <j_2\atop i_2-i_1<j_1-i_2} (\alpha(j_1-i_2))^{\delta/(1+\delta)}\\
&=& O\Big(\frac{(\log n)^{5r}}{nc_n^{2\delta/(1+\delta)}}\Big)\sum_{j=1}^n j(\alpha(j))^{\delta/(1+\delta)} \;=\; o(1).
\end{eqnarray*}
For the case $i_1=i_2$ we exemplarily consider the subcase $i_1=i_2<j_1<j_2$, other subcases are treated similarly. Note that $E[g(\xi_{i_1},\xi_{i_1},\xi_{j_1},\xi_{j_2}^*)]=0$ by the definition of $\varphi(\cdot)$, and $E[|g(\xi_{i_1},\xi_{i_1},\xi_{j_1},\xi_{j_2}^*)|^{1+\delta}]=O((\log n)^{\tilde r}c_n^{-2\delta})$ as before. Thus we obtain
\begin{eqnarray*}
\Bigg|\frac{1}{n^{3}}\sum_{i_1 <j_1 <j_2} E\Big[g(\xi_{i_1},\xi_{i_1},\xi_{j_1},\xi_{j_2})\Big]\Bigg|
&=& O((\log n)^{{\tilde r}/(1+\delta)})\frac{1}{n^3c_n^{2\delta/(1+\delta)}}\sum_{i_1 <j_1 <j_2} (\alpha(j_2-j_1))^{\delta/(1+\delta)}\\
&=& O\Big(\frac{(\log n)^{5r}}{nc_n^{2\delta/(1+\delta)}}\Big)\sum_{j=1}^n (\alpha(j))^{\delta/(1+\delta)} \;=\; o(1).
\end{eqnarray*}
Finally, the cases where more than two indices in $i_1,i_2,j_1,j_2$ are equal always lead to negligible terms by direct calculation. E.\,g.\ consider the term for $j_1=i_1\neq j_2=i_2$ in the sum (\ref{sum_u}). Applying assumption \ref{Am} its absolute value can straightforwardly be bounded by
$n^{-1}O((\log n)^{4r})(E[|\eps_1|])^2K^2(0)/c_n^2=o(1)$
by assumption \ref{Aneu}. The remaining terms are treated analogously.
 \hfill $\Box$

%\vspace{1cm}

%{\bf Acknowledgements.}

\end{appendix}

\end{document}